%% file: 0_main.tex
\documentclass[12pt,a4paper]{article}

\usepackage{array}
\usepackage{longtable}

\usepackage{amsmath}
\usepackage{amssymb}
\usepackage{amsthm}

\usepackage{color}

\usepackage{authblk}

\usepackage{graphicx}

\allowdisplaybreaks

\newtheoremstyle{plain_2}
  {\topsep}   % ABOVESPACE
  {\topsep}   % BELOWSPACE
  {\slshape}  % BODYFONT
  {0pt}       % INDENT (empty value is the same as 0pt)
  {\bfseries} % HEADFONT
  {.}         % HEADPUNCT
  {5pt} 		% HEADSPACE
  {}          % CUSTOM-HEAD-SPEC

\theoremstyle{definition}
\newtheorem{defi}{Definition}[section]
\newtheorem{rema}[defi]{Remark}
\newtheorem{exam}[defi]{Example}
\newtheorem{assu}[defi]{Assumption}

\newtheorem*{proo}{Proof}

\theoremstyle{plain_2}
\newtheorem{lemm}[defi]{Lemma}
\newtheorem{prop}[defi]{Proposition}
\newtheorem{theo}[defi]{Theorem}
\newtheorem{coro}[defi]{Corollary}

\numberwithin{equation}{section}

\DeclareMathOperator{\supp}{supp}
\DeclareMathOperator*{\Id}{Id}

\newcommand{\Div}{\mathrm{div}}

\newcommand{\trc}{\mathrm{tr}}
\newcommand{\Bil}{\mathrm{Bil}}

\newcommand{\Tt}{\vartheta}
\newcommand{\Tth}{\h{\vartheta}}
\newcommand{\ee}{\varepsilon}
\newcommand{\f}{\varphi}

\newcommand{\R}{\mathbb{R}}
\newcommand{\NN}{\mathbb{N}}

\newcommand{\Lt}{{L^2}}

\newcommand{\V}{{H_\Sigma^1}}
\newcommand{\Q}{{L^2}}

\newcommand{\pa}{\partial}
\newcommand{\G}{\nabla}
\newcommand{\Gh}{\h{\nabla}}

\newcommand{\Divh}{\h{\Div}}

\newcommand{\h}[1]{\widehat{#1}}
\newcommand{\ov}[1]{\overline{#1}}
\newcommand{\N}[1]{\left\| #1 \right\|}
\newcommand{\EN}[1]{\left| #1 \right|}
\newcommand{\DP}[2]{\left\langle #1 \right\rangle_{ #2 '\times #2 }}

\newcommand{\Of}{\int_{\Omega_f}}
\newcommand{\Os}{\int_{\Omega_s}}

\newcommand{\Ost}{\int_{\Omega_s(t)}}

\newcommand{\dxh}{\,d\h{x}}
\newcommand{\dx}{\,dx}
\newcommand{\dt}{\,dt}

\usepackage{csquotes}
\usepackage[backend=biber,style=numeric,uniquename=init,giveninits=true,isbn=false,url=false,doi=false,eprint=false]{biblatex}
\usepackage[colorlinks,citecolor=black,linkcolor=black,urlcolor=black,filecolor=black,menucolor=black]{hyperref}

\usepackage[left=2cm, right=2cm, top=2.5cm, bottom=3cm, footskip=2cm]{geometry}

\title{Derivation and analysis of a Stokes-transport system in evolving vessels modeling thermoregulation in human skin}
\author[1,2]{Kilian Hacker}
\author[1]{Maria Neuss-Radu}
\affil[1]{Department of Mathematics, Friedrich-Alexander-Universit\"{a}t Erlangen-N\"{u}rnberg, Cauerstr.~11, 91058 Erlangen, Germany.}
\affil[2]{Corresponding author. E-mail address: \href{mailto:kilian.hacker@fau.de}{kilian.hacker@fau.de}}
\date{February 9, 2026}

\begin{document}

\maketitle

\begin{abstract}
    We consider a Stokes flow coupled with advective-diffusive transport in an evolving domain with boundary conditions allowing for inflow and outflow. The evolution of the domain is induced by the transport process, leading to a fully coupled problem. Our aim is to model the thermal control of blood flow in human skin. To this end, the model takes into account the temperature-dependent production of biochemical substances, the subsequent dilation and constriction of blood vessels, and the resulting changes in convective heat transfer. We prove existence and uniqueness of weak solutions using a fixed point method that allows us to treat the nonlinear coupling.
\end{abstract}

\textbf{Keywords.} Stokes equations, heat transport, free boundary problem, coupled problem, vasodilation\par\bigskip
\textbf{MSC2020 subject classification.} 35R35, 35K57, 76D07, 80A19, 92C30
%\par\bigskip\textbf{Abbreviated title.} A Stokes-transport system for thermoregulation

\input{1_introduction}

\input{2_model}

\input{3_result}

\input{4_proof}

\section*{Conclusion and Outlook}
In this paper, a new mathematical model for thermoregulation in human skin was derived and analyzed. The key components are blood flow dynamics and heat transfer, as well as the biochemical processes in the arterial wall that lead to vasodilation and vasoconstriction. Moreover, the full coupling of these subsystems was taken into account. Consequently, numerical simulations based on the proposed model can provide a more comprehensive understanding of thermoregulation in human skin. \par
At the same time, the present model can serve as a starting point for future research on the process of thermoregulation. As a first generalization, the model could be extended by incorporating the elastic behavior of the arterial wall and the surrounding tissue, which would lead to a fluid-structure interaction problem. We note that in the context of atherosclerosis models, fluid-structure interaction in connection with transport and growth processes has been considered, for example, in \cite{yang_mathematical_2016} and \cite{abels_fluid-structure_2023}. Furthermore, modeling the effect of the wall shear stress on the NO production in the endothelium would also be a valuable extension of the model, see also \cite{kudryashov_numerical_2008, wang_numerical_2019, fok_shear_2023, marino_unraveling_2024}. In both cases, the resulting models pose major challenges for analysis and simulation.

\section*{Acknowledgements}
The authors thank Willi Jäger for fruitful discussions on the topic.

\section*{Research funding}
The research of K.H.~was funded by the Deutsche Forschungsgemeinschaft (DFG, German Research Foundation) -- GRK2339 -- Projektnummer 321821685.

\printbibliography
\end{document}

%% file: 1_introduction.tex
%Introduction
\section{Introduction} \label{sec:introduction}

The aim of human thermoregulation is to maintain a constant core body temperature. The control of skin blood flow plays a key role in this process: Vasodilation in human skin increases convective heat transfer from the body core to the surface, while vasoconstriction reduces heat loss to the environment \cite{charkoudian_mechanisms_2010}. 
In this paper, we describe these control mechanisms by means of a mathematical model in the form of partial and ordinary differential equations in evolving domains, which takes into account the interplay between blood flow dynamics, heat transfer as well as biochemical processes.

%\paragraph{Mechnisms of thermoregulation}
Arterioles and arteries in the skin are lined with an inner layer of endothelial cells and several layers of smooth muscle cells. Under body heating or cooling, the smooth muscle cells are able to dilate and constrict the vessels, thereby regulating skin blood flow levels. In our investigations, we focus on the cases where the human body is only locally exposed to heating or cooling, for example on the skin of the forearm, where many experimental studies have been conducted \cite{terjung_cutaneous_2014}. As pointed out in \cite[p.~55]{terjung_cutaneous_2014}, \textquotedblleft our current understanding of the response in skin blood flow to local heating now includes complex interactions between neural components and locally produced chemical messengers\textquotedblright. Hereby, the molecule nitric oxide (NO) plays a key role in many of these pathways. The production of NO takes place in endothelial cells, catalyzed by the enzyme NO synthase (NOS). NOS is activated by several neurotransmitters during local heating. An increased NO concentration in the endothelium induces relaxation of the smooth muscle cells, leading to vessel dilation. For comprehensive overviews of these physiological findings, see, for example,  \cite{charkoudian_mechanisms_2010}, \cite{johnson_local_2010}, and \cite{terjung_cutaneous_2014}.

%\paragraph{Summary of the mathematical model}
We develop a mathematical model for thermoregulation in human skin describing the processes in one single artery or arteriole and in the surrounding cutaneous tissue. The system of equations comprises a quasi-stationary Stokes system for the blood flow, an advection-diffusion equation for the evolution of temperature, and an ordinary differential equation (ODE) for the concentration of NO. We formulate the system in a tube-shaped fluid domain and a surrounding solid domain, separated by an evolving interface, which represents the vessel wall. The equations are fully and nonlinearly coupled: The Stokes system provides the convective transport velocity for the advection-diffusion equation, the temperature influences the production of NO, and the concentration of NO determines the evolution of the fluid/solid interface. Consequently, the deformation of the whole domain is an unknown of the model and contributes to the nonlinear coupling, see also Figure \ref{fig:coupling} (left) for a schematic representation of the couplings in the model. In addition, we pay special attention to adequate boundary and transmission conditions. In particular, the Stokes flow satisfies a no-slip condition at the fluid/solid interface and a normal stress condition at the ends of the fluid domain, allowing for inflow and outflow. A Danckwerts boundary condition for the normal heat flux accounts for the heat transported by the inflowing blood, while a Dirichlet boundary condition at the solid part models some external heating. Moreover, at the fluid/solid interface, the normal heat flux is continuous and proportional to the jump in temperature.

%\paragraph{Summary of methods and contributions}
We prove global-in-time existence and uniqueness of a solution to the fully coupled system of equations by means of Schaefer's fixed point theorem. As the domain is evolving in time, we first transform the weak formulation of the problem to a reference domain. This leads to a nonlinear coupling through coefficients in all equations, involving the unknown deformation gradient. For the Stokes system, we get additional source terms by subtracting the boundary data from the unknowns. Existence and uniqueness then follows from a generic result for saddle point problems. For the advection-diffusion equation, an additional challenge arises from the coupling of the temperatures in the fluid and the solid part by a transmission condition at the fluid/solid interface. Here, we adapt the standard Galerkin method so that it captures both the two unknowns and the coefficients. Finally, we provide all estimates necessary to prove the continuity of the solution operator and the uniqueness of the fully coupled problem.

%\paragraph{Literature related by the modeling}
To the best of our knowledge, the model for thermoregulation presented in this paper is the first of its kind to include all the physiological processes mentioned above. In fact, there are only a few comprehensive mathematical models in the literature for the local thermal control of blood flow. We mention here the contributions in \cite{kudryashov_numerical_2008, wang_numerical_2019}. The model developed in \cite{kudryashov_numerical_2008} considers a single axially symmetric and viscoelastic artery. The radius of the artery wall is assumed to be constant over the vessel segment, so the blood flow can be described by a quasi-stationary (generalized) Poiseuille law. The equation of motion of the vessel wall incorporates the active force exerted by smooth muscle cells. The smooth muscle activity is determined by biochemical processes including the wall shear stress mediated production of NO in the endothelial layer, see also \cite{fok_shear_2023, marino_unraveling_2024} for other numerical investigations on this topic. In \cite{wang_numerical_2019}, the model from \cite{kudryashov_numerical_2008} is placed in the context of thermoregulation by assuming different blood flow rates at different temperatures. A significant limitation of this modeling approach is its inability to provide feedback from changes in blood flow to temperature, despite the fact that this is an essential aspect of thermoregulation. The limitations in the above models were partly overcome in \cite{schorten_analysis_2013}, where an advection-diffusion equation for the temperature of blood and tissue is coupled to the blood flow and the movement of the vessel wall. However, although the vessel radius may change along the vessel axis, the blood flow is described by a potential flow instead of, e.g., a Stokes flow. Furthermore, this model does not include the biochemical component, involving, e.g., the NO dynamics.

%\paragraph{Literature related by the analysis}
The major analytical challenge of this work is to deal with the coupling of the equations. Related problems arise, e.g., in precipitation and dissolution models: The microscopic model under consideration in \cite{gahn_rigorous_2025} includes a Stokes flow, an advection-diffusion problem, and an ODE for the evolution of the domain, see also the previous works \cite{gahn_homogenization_2021, gahn_homogenization_2023, wiedemann_homogenisation_2023, wiedemann_homogenisation_2024}. The transformed problem in the reference domain is fully coupled via coefficients including the unknown deformation gradient. However, the transport process in \cite{gahn_rigorous_2025} only takes place in the fluid part. A model for fluid flow with solute transport in both fluid and solid is formulated in \cite{quarteroni_mathematical_2002}. At the fluid/solid interface, the model states that the normal solute fluxes are equal and, moreover, proportional to the difference of concentrations in the fluid and solid part, like in our model. The major difference to the present work is that the domain does not change. In the work \cite{mabuza_modeling_2016}, the authors consider a reactive transport problem in a fluid and an adsorption-desorption problem on the boundary of the fluid domain, which is evolving in time. However, the transport problem has no effect on the fluid flow or the evolution of the domain, as the corresponding fluid-structure interaction problem is solved separately, see \cite{muha_existence_2013}. With our model, we are now able to treat all these couplings at the same time.

%\paragraph{Structure of the paper}
The paper is organized as follows. The mathematical model is formulated in Sections \ref{subsec:domain} to \ref{subsec:ode} based on insights regarding the physiological processes. In Section \ref{subsec:weak_formulation}, we derive the weak formulation in the reference domain and outline all assumptions on the data. The main existence and uniqueness theorem is stated and proven in Section \ref{sec:result}.

\paragraph{Notation}
Given a matrix $M\in\R^{n\times n}$, we write $M^T$ for its transposed, and $M^{-T}$ short hand for $(M^{-1})^T$ if $M$ is invertible. The identity matrix of size 3 is denoted by $I$, the Euclidean norm of $v\in\R^n$ by $|v|$, and the Frobenius inner product of $M,N\in\R^{n\times n}$ by $M:N := \trc(M^TN)$. \par
For any real-valued function $f$, we set $f^-=\min\{f,0\}$ and $f^+=\max\{f,0\}$. For a differentiable vector field $v:\R^n\to\R^n$, let $Dv$ be its Jacobian, i.e.~$(Dv)_{ij} = \pa_{x_j}v_i$, and define $\G v := (Dv)^T$. Moreover, for a differentiable function $M:\R^n\to\R^{n\times n}$, the divergence $\Div M$ is a column vector with the components $(\Div M)_i := \sum_{j=1}^{n}\pa_{x_j}M_{ji}$, so that it holds $\Delta v = \Div(\G v)$. \par
In Section \ref{sec:model}, the mathematical model will be formulated in the evolving domain, i.e.~in the Eulerian framework. We denote the spatial variable associated with the Eulerian framework by $\h{x}$ and add a hat to all quantities depending on $\h{x}$ and to all derivatives with respect to $\h{x}$, e.g., $\h{f}$ or $\Gh$. Moreover, we write:
\begin{equation*}
    \textstyle\h{e}(\h{v}):=\frac{1}{2}\big(\Gh\h{v}+(\Gh\h{v})^T\big)\qquad\text{and}\qquad e_M(v):=\frac{1}{2}\big(M^{-T}\G v + (M^{-T}\G v)^T\big).
\end{equation*}\par
If $\Omega\subset\R^3$ is a bounded open Lipschitz domain and $\Gamma\subset\pa\Omega$ with $\int_\Gamma d\sigma >0$, we denote by $H^1_\Gamma(\Omega)$ the subspace of $H^1(\Omega)$ of functions with zero trace on $\Gamma$. Finally, for functions of $t$ and $x_1$ or $x$, and for $T,L>0$, we introduce the function spaces
\begin{align*}
    C_1^2([0,T]\times[0,L]) &:= \{f:[0,T]\times[0,L]\to\R: f,\pa_t f,\pa_{x_1}f,\pa_{x_1}\pa_t f, \pa_{x_1}^2 f\in C([0,T]\times[0,L])\}, \\
    C_1^2([0,T]\times\ov{\Omega}) &:= \{f:[0,T]\times\ov{\Omega}\to\R: \pa_t^{\beta_1}\pa_{x_1}^{\beta_2}\pa_{x_2}^{\beta_3}\pa_{x_3}^{\beta_4}f\in C([0,T]\times\ov{\Omega})\text{ for all  }\beta\in\NN_0^4 \\
    &\hspace{7cm}\text{ with }0\leq |\beta|\leq 2\text{ and }\beta_1\leq 1 \}.
\end{align*}

%% file: 2_model.tex
%Model
\section{Mathematical model} \label{sec:model}

\begin{figure}[tbh]
    \centering
    \includegraphics[width=0.45\linewidth]{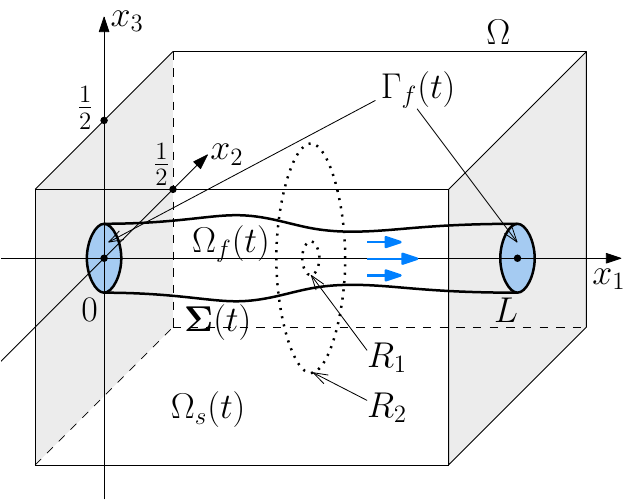} \hfill
    \includegraphics[width=0.45\linewidth]{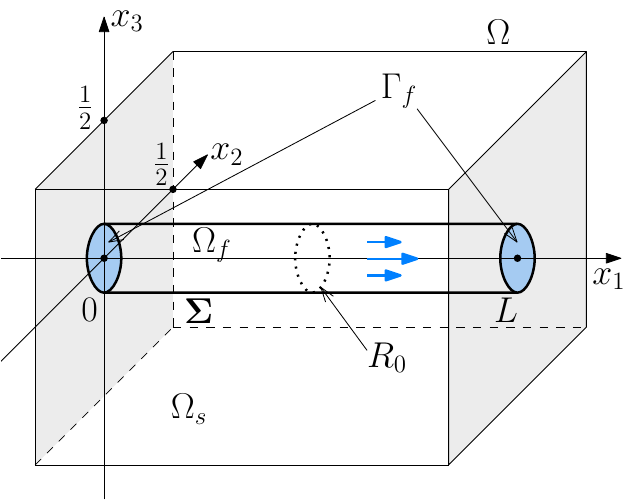}
    \caption{The evolving domain (left) and the reference domain (right), both divided into a fluid and a solid part.}
    \label{fig:domain}
\end{figure}

We aim to describe the behavior of a single artery or arteriole in human skin. Therefore, we consider fluid flow, heat transport, and the concentration of a chemical substance in an evolving domain, which is shown in Figure \ref{fig:domain}. In Sections \ref{subsec:domain} to \ref{subsec:ode}, the mathematical model in the evolving domain is presented. The weak formulation of the problem in the reference domain is derived in Section \ref{subsec:weak_formulation}.

\subsection{The evolving domain} \label{subsec:domain}

We model one single blood vessel and the surrounding tissue by a tube-shaped fluid domain and a cuboid solid domain, both changing in time. The fluid/solid interface represents the vessel wall. We allow this interface to move, which corresponds to the dilation and constriction of the blood vessel. This causes a deformation of the fluid and solid domain. \par
Fix a minimal and a maximal radius $0<R_1<R_2<\frac{1}{2}$, a reference radius $R_0\in[R_1,R_2]$, some $\delta>0$ with $R_2+3\delta<\frac{1}{2}$ and $R_1-3\delta>0$, and a length $L>0$. In this section, we employ the notation $x=(x_1,\bar{x})\in\R\times\R^2$. \par
Let the reference domain be $\Omega := \left(0,L\right)\times\left(-\frac{1}{2},\frac{1}{2}\right)^2$ and define the fluid and solid subdomains by
\begin{equation*}
    \Omega_f := \left\{(x_1,\bar{x})\in\Omega: |\bar{x}|<R_0\right\}\qquad\text{and}\qquad
    \Omega_s:=\Omega\setminus\overline{\Omega_f},
\end{equation*}
which are separated by the interface
\begin{equation*}
    \Sigma:=\left\{(x_1,\bar{x})\in\ov{\Omega}: |\bar{x}|=R_0\right\}.
\end{equation*}
In accordance with the mixed boundary conditions formulated below, we divide the boundary $\pa\Omega$ into three parts, namely the inflow/outflow part
\begin{equation*}
    \Gamma_f:=\left\{(x_1,\bar{x})\in\ov{\Omega}: x_1\in\{0,L\}\text{ and }|\bar{x}|\leq R_0\right\},
\end{equation*}
the Dirichlet part $\Gamma_D:=[0,L]\times\left[-\frac{1}{2},\frac{1}{2}\right]\times\left\{\frac{1}{2}\right\}$ and the remaining part $\Gamma_N=\ov{\pa\Omega\setminus(\Gamma_f\cup\Gamma_D)}$. \par
Now, we introduce the evolving subdomains by means of the deformation $S$ as follows. The radius $R:[0,T]\times[0,L]\to[R_1,R_2]$ of the evolving interface results from the solution of the ODE by the relation \eqref{eq:R_from_c}. The deformation for the whole domain is then defined by
\begin{equation} \label{eq:deformation_def}
    S:[0,T]\times\ov{\Omega}\to\R^3,~S(t,x)=\begin{cases}\begin{pmatrix}
	   x_1 \\
	   \rho(R(t,x_1),\EN{\bar{x}})\frac{\bar{x}}{\EN{\bar{x}}}
    \end{pmatrix}&\qquad\text{if }|\bar{x}|\neq 0, \\
    x &\qquad\text{if }|\bar{x}|=0,
\end{cases}
\end{equation}
with a suitable function $\rho:[R_1,R_2]\times[0,\infty)\to[0,\infty)$, see Assumption \ref{assu:data} (A1). Finally, we define the evolving subdomains by $\Omega_f(t) := S(t,\Omega_f)$ and $\Omega_s(t):= S(t,\Omega_s)$, the evolving interface by $\Sigma(t):= S(t,\Sigma)$, and the evolving boundaries by $\Gamma_f(t):= S(t,\Gamma_f)$, $\Gamma_N(t):= S(t,\Gamma_N)$, and $\Gamma_D(t):= S(t,\Gamma_D)$, see also Remark \ref{rema:S_bijectivity}.

\subsection{Stokes system modeling blood flow} \label{subsec:stokes}

Since the model concerns only medium vessels and no capillaries, blood can approximately be considered as a Newtonian fluid, see, e.g., \cite{quarteroni_computational_2000} for a discussion of blood rheology modeling. We neglect the pulsating behavior of blood flow and assume that the non-periodic variation in blood pressure and vessel diameter are much slower than the flow itself. Thus, a quasi-stationary Stokes flow is an acceptable approximation. At the outer boundary of the fluid domain, we allow for inflow and outflow under normal stress conditions, since the fluid domain represents a single vessel within a larger network. See also \cite{fabricius_stokes_2019} for a mathematical discussion of appropriate boundary conditions in this case. \par
We consider the following quasi-stationary Stokes problem: Find $\h{v}_f:\bigcup_{t\in[0,T]}\{t\}\times\Omega_f(t)\to\R^3$ and $\h{p}:\bigcup_{t\in[0,T]}\{t\}\times\Omega_f(t)\to\R$ solving
\begin{equation} \label{eq:stokes_eq}
    \left\{\enspace\begin{aligned}
        &\Divh(-\h{p}\,I + 2\mu\,\h{e}(\h{v}_f))=0 & \qquad\text{in }\Omega_f(t),&\\
        &\Divh\h{v}_f=0 & \text{in }\Omega_f(t),&\\
    \end{aligned}\right.
\end{equation}
with the normal stress boundary condition at the inflow/outflow boundary:
\begin{equation} \label{eq:stokes_bc_inout}
    (-\h{p}\,I + 2\mu\,\h{e}(\h{v}_f))\h{n}_{t}=-\h{f}_b\h{n}_{t}\qquad\text{on }\Gamma_f(t),
\end{equation}
and the no-slip condition at the fluid/solid interface:
\begin{equation} \label{eq:stokes_bc_interface}
    \h{v}_f(t,\cdot)=\pa_t S\left(t,S^{-1}(t,\cdot)\right)\qquad\text{on }\Sigma(t),
\end{equation}
for all $t\in[0,T]$. \par
We denote by $\h{n}_{t}$ the outer normal to $\partial\Omega_f(t)$, by $\mu$ the dynamic viscosity of the fluid, and by $\h{f}_b:\bigcup_{t\in[0,T]}\{t\}\times\Gamma_f(t)\to\R$ the boundary data for the normal stress.

\begin{rema} \label{rema:inflow_outflow}
    If $\h{f}_b$ represents a pressure at the boundary and the pressure drop between $x_1=0$ and $x_1=L$ is sufficiently high, we expect an inflow at $x_1=0$ and an outflow at $x_1=L$. However, this is not clear a priori and is not necessary for the well-posedness of the model. For example, a fast dilation of $\Omega_f(t)$ could cause an inflow both at $x_1=0$ and $x_1=L$.
\end{rema}

\subsection{Advection-diffusion equation modeling temperature evolution} \label{subsec:adv_diff}

We describe the heat transport in blood and tissue by a system of advection-diffusion equations. A Dirichlet boundary condition at the solid part accounts for some external heating of the tissue, and a Danckwerts boundary condition at the fluid part incorporates the heat transported by the inflowing blood. At the fluid/solid interface, we impose a transmission condition that allows for a jump in temperature, since the interface represents the vessel wall consisting of several layers of endothelial and smooth muscle cells. In addition, the advection-diffusion equations are coupled to the Stokes flow, so we are able to investigate the influence of vasodilation on heat transport. \par
The evolution of temperature is subject to the following transport problem: Find $\Tth_f:\bigcup_{t\in(0,T)}\{t\}\times\Omega_f(t)\to\R$ and $\Tth_s:\bigcup_{t\in(0,T)}\{t\}\times\Omega_s(t)\to\R$ satisfying
\begin{equation} \label{eq:addiff_eq}
    \left\{\enspace\begin{aligned}
        \textstyle&\pa_t\Tth_f + \Divh(\h{v}_f\Tth_f-K_f\Gh\Tth_f)=0 & \qquad\text{in }\bigcup_{t\in(0,T)}\{t\}\times\Omega_f(t),&\\
        \textstyle&\pa_t\Tth_s + \Divh(\h{v}_s\Tth_s-K_s\Gh\Tth_s)=0 & \qquad\text{in }\bigcup_{t\in(0,T)}\{t\}\times\Omega_s(t),&\\
    \end{aligned}\right.
\end{equation}
where
\begin{equation}\label{eq:v_s_def}
\h{v}_s(t,\cdot):=\pa_t S\left(t,S^{-1}(t,\cdot)\right)\qquad\text{in }\Omega_s(t),
\end{equation}
with the transmission condition of imperfect heat transfer
\begin{equation} \label{eq:addiff_tr}
    -K_f\Gh\Tth_f\cdot\h{n}_{t}=-K_s\Gh\Tth_s\cdot\h{n}_{t}=\alpha(\Tth_f-\Tth_s)\qquad\text{on }\bigcup_{t\in(0,T)}\{t\}\times\Sigma(t),
\end{equation}
the Danckwerts boundary condition (see \cite{danckwerts_continuous_1953})
\begin{equation} \label{eq:addiff_bc_danckwerts}
    \Tth_f\h{v}_f\cdot\h{n}_{t}-K_f\Gh\Tth_f\cdot\h{n}_{t} = \Tth_f(\h{v}_f\cdot\h{n}_{t})^+ + \h{f}_{in}(\h{v}_f\cdot\h{n}_{t})^- \qquad\text{on }\bigcup_{t\in(0,T)}\{t\}\times\Gamma_f(t),\\
\end{equation}
the mixed boundary conditions
\begin{align}
    -K_s\Gh\Tth_s\cdot\h{n}_{t}&=0 & &\text{on }\bigcup_{t\in(0,T)}\{t\}\times\Gamma_{N}(t), \label{eq:addiff_bc_neumann}\\
    \Tth_s&=\h{\Tt}_{ext} & &\text{on }(0,T)\times\Gamma_{D}, \label{eq:addiff_bc_dirichlet}
\end{align}
and the initial conditions $\Tth_f|_{t=0}=\Tth_{f}^0$ in $\Omega_f(0)$ and $\Tth_s|_{t=0}=\Tth_{s}^0$ in $\Omega_s(0)$. \par
The matrices $K_f,K_s\in\R^{3\times 3}$ represent the thermal conductivities of the fluid and the solid, $\h{f}_{in}:\bigcup_{t\in(0,T)}\{t\}\times\Gamma_f(t)\to\R$ is the boundary data for the heat flux at $\{x\in\Gamma_f(t):\h{v}_f(t)\cdot\h{n}_t\leq 0\}$, $\h{\Tt}_{ext}:[0,T]\times\Gamma_D\to\R$ is the temperature of the external heating at $\Gamma_D$, $\alpha>0$ is the heat transfer coefficient on $\Sigma(t)$, and $\h{n}_{t}$ denotes the outer normal to $\partial\Omega_f(t)$ (on $\Sigma(t)$) and to $\partial\Omega$.

\begin{rema} \label{rema:addiff_bc_danckwerts}
The boundary condition \eqref{eq:addiff_bc_danckwerts} reduces to $-K_f\Gh\Tth_f\cdot\h{n}_{t}=0$ on the parts of $\Gamma_f(t)$ where $\h{v}_f\cdot\h{n}_{t}\geq 0$. However, \eqref{eq:addiff_bc_danckwerts} will prove useful in the weak formulation. 
\end{rema}

\subsection{ODE modeling chemical reactions} \label{subsec:ode}

During local heating, the enzyme NOS is activated in the endothelium and produces NO \cite{terjung_cutaneous_2014}. We model this processes by an ODE for the NO concentration at the fluid/solid interface, with a production term depending on the average $\h{\mathcal{T}}(\h{\Tt}_s)$ of the temperature in the solid domain. Thus, the data for the ODE depends on $t$, and additionally on $x_1$, which models the distribution of NOS within the
endothelium. As a consequence, we search for a solution $c:[0,T]\times[0,L]\to\R$ of the ODE
\begin{equation} \label{eq:ode}
    \pa_t c(t,x_1) = -k\,c(t,x_1) + G(x_1,\h{\mathcal{T}}(\h{\Tt}_s)(t))\hspace{1cm}\text{for all }(t,x_1)\in[0,T]\times[0,L],
\end{equation} 
with the initial condition
\begin{equation} \label{eq:initial_concentration}
    c(0,x_1)=c^0(x_1)\qquad\text{for all }x_1\in[0,L].
\end{equation}
Here, $k>0$ is the degradation rate constant, and $G:[0,L]\times\R\to\R$ is the production rate function, see Assumption \ref{assu:data} (A2).
 
We define the averaging operator $\h{\mathcal{T}}$ for some fixed $0<\gamma<T$ in two steps: First, we take a weighted average of $\h{\Tt}_s$ over $\Omega_s(t)$ and extend it in time to $(-\gamma,0]$ by the weighted average of the  initial data:
\begin{equation} \label{eq:space_averaged_temp}
    \begin{aligned}
        &\h{\mathcal{T}}_1: \Lt(0,T;\Lt(\Omega_s(t)))\to\Lt(-\gamma,T),\\
        &\qquad\h{\mathcal{T}}_1(\h{\Tt}_s)(t)=\begin{cases}
            \frac{1}{|\Omega_s|}\int_{\Omega_s(0)}\frac{\h{\Tt}_s^0(\h{x})}{\det(DS(0,S^{-1}(0,\h{x})))}\dxh &\text{if }t\in(-\gamma,0], \\
            \frac{1}{|\Omega_s|}\Ost\frac{\h{\Tt}_s(t,\h{x})}{\det(DS(t,S^{-1}(t,\h{x})))}\dxh &\text{if }t\in(0,T).
        \end{cases}
    \end{aligned}
\end{equation}
Second, we take a weighted average of $\h{\mathcal{T}}_1(\h{\Tt}_s)$ over the time interval $(t-\gamma,t]$:
\begin{equation} \label{eq:time_averaged_temp}
    \begin{aligned}
        &\h{\mathcal{T}}: \Lt(0,T;\Lt(\Omega_s(t)))\to C([0,T]),\\
        &\qquad\h{\mathcal{T}}(\h{\Tt}_s)(t) = \int_{\R}K_\gamma(t-s)\,\h{\mathcal{T}}_1(\h{\Tt}_s)(s)\,ds = (K_\gamma\star\h{\mathcal{T}}_1(\h{\Tt}_s))(t),
    \end{aligned}
\end{equation}
where $K_\gamma\in C_c(\R)$ is a non-negative function with $\supp(K_\gamma)\subset[0,\gamma)$ and $\int_{\R}K_\gamma(t)\,dt = 1$. 
The averaging \eqref{eq:space_averaged_temp} is motivated by the averaged temperature registered by sensory nerves within the tissue leading to a nonlocal dependence of the NO-production on temperature. Modeling approaches that use nonlocal dependencies can be also found, for example, in \cite{jager_analysis_2009} for biochemical processes in living cells, and \cite{gahn_rigorous_2025}, for precipitation and dissolution processes in the geosciences. \par
Finally, we model the dependence of the vessel radius on the NO concentration using a function $H:\R\to[R_1,R_2]$ that incorporates the dynamics of smooth muscle cells in a phenomenological manner. 
In particular, we use $H$ to set a lower limit $R_1$ and an upper limit $R_2$ for the vessel radius. This reflects the fact that the diameter of the blood vessel should neither become arbitrarily large nor shrink to zero. Thus, the radius of the evolving interface is given by 
\begin{equation} \label{eq:R_from_c}
    R:[0,T]\times[0,L]\to[R_1,R_2],~R(t,x_1)=H(c(t,x_1)),
\end{equation}
see also Assumption \ref{assu:data} (A3). Note that, in general, $R(0,x_1)\neq R_0$.

\subsection{Weak formulation in the reference domain} \label{subsec:weak_formulation}

In this section, we give the main steps for the derivation of the weak formulation of the mathematical model, which is presented in Definition \ref{defi:weak_formulation}. First, we derive the weak formulation for the Stokes system and the advection-diffusion equation in the evolving domain. Then, using the deformation $S$, we transform the model to the reference domain. In doing so, the general notation is the following: The inverse of $S(t,\cdot)$ is denoted by $S^{-1}(t,\cdot)$ for any $t\in[0,T]$, the spatial variables $\h{x}\in M(t)$ and $x\in M$ are related by $\h{x}=S(t,x)$ whenever $M(t)=S(t,M)$, the quantities $\h{f}$ and $f$ are related by $f(t,x)=\h{f}(t,S(t,x))$ for all $x\in M$, and the derivatives $\Divh$ and $\Gh$ are with respect to $\h{x}$, while $\Div$ and $\G$ are with respect to $x$. Moreover, we use the notation \eqref{eq:coefficients} for the coefficients. \par
To formulate the ODE in the reference domain, it suffices to define the operator $\mathcal{T}$ by
\begin{equation*}
    \mathcal{T}: L^2(0,T;L^2(\Omega_s))\to C([0,T]),\qquad\mathcal{T}(\Tt_s) := \h{\mathcal{T}}(\h{\Tt}_s) = \h{\mathcal{T}}(\Tt_s\circ S^{-1}).
\end{equation*}
We then obtain the representation \eqref{eq:T_representation} by transforming the integrals in the definition \eqref{eq:space_averaged_temp} of $\h{\mathcal{T}}_1$ to the reference domain. \par
Concerning the momentum equation of the Stokes system \eqref{eq:stokes_eq}, let us first define
\begin{equation} \label{eq:v_b_def}
\h{v}_b(t,\cdot):=\pa_t S\left(t,S^{-1}(t,\cdot)\right)\qquad\text{in }\Omega_f(t),
\end{equation}
which extends the data for \eqref{eq:stokes_bc_interface} from $\Sigma(t)$ to $\Omega_f(t)$, and choose $\h{f}_b:\bigcup_{t\in[0,T]}\{t\}\times\Omega_f(t)\to\R$ so that \eqref{eq:stokes_bc_inout} holds in the trace sense. Now, we make the substitutions $\h{v}_f=\h{w}+\h{v}_b$ and $\h{p}=\h{q}+\h{f}_b$, that is, we subtract the boundary data from the unknowns. Multiplying the momentum equation by a test function $\h{\phi}:\Omega_f(t)\to\R^3$ with $\h{\phi}=0$ on $\Sigma(t)$, integrating over $\Omega_f(t)$, integrating by parts and using the boundary condition \eqref{eq:stokes_bc_inout}, we obtain:
\begin{equation*}
    \int_{\Omega_f(t)}(-\h{p}\,I + 2\mu\,\h{e}(\h{v}_f)):\Gh\h{\phi}\dxh = -\int_{\Gamma_f(t)}\h{f}_b\h{\phi}\cdot\h{n}_{t}\,d\h{\sigma}.
\end{equation*}
By the above substitutions and again by integration by parts, it follows:
\begin{equation*}
    \int_{\Omega_f(t)}(-\h{q}\,I + 2\mu\,\h{e}(\h{w})):\Gh\h{\phi}\dxh = -\int_{\Omega_f(t)}\Gh\h{f}_b\cdot\h{\phi} + 2\mu\,\h{e}(\h{v}_b):\Gh\h{\phi}\dxh,
\end{equation*}
and with the same substitutions, we get $\h{\Div}(\h{w})=-\h{\Div}(\h{v}_b)$. Finally, a change of variables yields \eqref{eq:stokes_weak_ref}. Note that $|J|=J$ by \eqref{eq:S_estimate_J} and $A:\G w = \Div(Aw)$ by the product rule and the Piola identity \eqref{eq:piola_identity}. \par
In order to derive the weak formulation for the advection-diffusion equation, we multiply \eqref{eq:addiff_eq} by test functions $\h{\f}_f:\Omega_f(t)\to\R$ resp. $\h{\f}_s:\Omega_s(t)\to\R$ with $\f_s=0$ on $\Gamma_D$, integrate over $\Omega_f(t)$ resp. $\Omega_s(t)$, add both equations, integrate by parts and use the conditions \eqref{eq:addiff_tr}, \eqref{eq:addiff_bc_danckwerts}, \eqref{eq:addiff_bc_neumann}, and \eqref{eq:addiff_bc_dirichlet} to obtain
\begin{align*}
    &\int_{\Omega_f(t)}\pa_t\h{\Tt}_f\h{\f}_f\dxh + \int_{\Omega_f(t)}(K_f\Gh\h{\Tt}_f - \h{v}_f\h{\Tt}_f)\cdot\Gh\h{\f}_f\dxh \\
    &+ \int_{\Omega_s(t)}\pa_t\h{\Tt}_s\h{\f}_s\dxh + \int_{\Omega_s(t)}(K_s\Gh\h{\Tt}_s - \h{v}_s\h{\Tt}_s)\cdot\Gh\h{\f}_s\dxh \\
    &+ \int_{\Sigma(t)}(\h{\Tt}_f\h{\f}_f - \h{\Tt}_s\h{\f}_s)\h{v}_s\cdot\h{n}_t\,d\h{\sigma} + \int_{\Sigma(t)}\alpha(\h{\Tt}_f - \h{\Tt}_s)(\h{\f}_f - \h{\f}_s)\,d\h{\sigma} \\
    &= -\int_{\Gamma_f(t)}(\h{\Tt}_f(\h{v}_f\cdot\h{n}_t)^+ + \h{f}_{in}(\h{v}_f\cdot\h{n}_t)^-)\h{\f}_f\,d\h{\sigma}.
\end{align*}
Note that by definitions \eqref{eq:v_s_def}, \eqref{eq:v_b_def}  and by \eqref{eq:stokes_bc_interface}, $\h{v}_f = \h{v}_b = \h{v}_s$ on $\Sigma(t)$. We write $\h{\pa_t\Tt_{f/s}}$ for the time derivatives in the reference domain with coordinates in the evolving domain, that is,
\begin{equation*}
    \h{\pa_t\Tt_{f/s}}(t,S(t,x)) := \pa_t\Tt_{f/s}(t,x),
\end{equation*}
so the chain rule yields the identities
\begin{align*}
    &\begin{aligned}
        \h{\pa_t\Tt_{f/s}}(t,S(t,x)) &= \frac{d}{dt}(\Tt_{f/s}(\cdot,x))(t) = \frac{d}{dt}(\h{\Tt}_{f/s}(\cdot,S(\cdot,x)))(t) \\
        &= \pa_t\h{\Tt}_{f/s}(t,S(t,x)) + \h{D}\h{\Tt}_{f/s}(t,S(t,x))\pa_t S(t,x)
    \end{aligned}\\
    \Longleftrightarrow~&\pa_t\h{\Tt}_{f/s}(t,\h{x}) = \h{\pa_t\Tt_{f/s}}(t,\h{x}) - \h{v}_{b/s}(t,\h{x})\cdot\Gh\h{\Tt}_{f/s}(t,\h{x}),
\end{align*}
where `$f/s$' indicates an equation for both `$f$' and `$s$'. Using integration by parts and the facts that $\h{v}_b\cdot\h{n}_t=0$ on $\Gamma_f(t)$ and $\h{v}_s\cdot\h{n}_t=0$ on $\Gamma_D\cup\Gamma_N(t)$, the advection-diffusion equation takes the form
\begin{align*}
    &\int_{\Omega_f(t)}\h{\pa_t\Tt_f}\h{\f}_f + \h{\Div}(\h{v}_b)\h{\Tt}_f\h{\f}_f\dxh + \int_{\Omega_f(t)}(K_f\Gh\h{\Tt}_f - (\h{v}_f-\h{v}_b)\h{\Tt}_f)\cdot\Gh\h{\f}_f\dxh \\
    &+ \int_{\Omega_s(t)}\h{\pa_t\Tt_s}\h{\f}_s + \h{\Div}(\h{v}_s)\h{\Tt}_s\h{\f}_s\dxh + \int_{\Omega_s(t)}K_s\Gh\h{\Tt}_s\cdot\Gh\h{\f}_s\dxh \\
    &+ \int_{\Sigma(t)}\alpha(\h{\Tt}_f - \h{\Tt}_s)(\h{\f}_f - \h{\f}_s)\,d\h{\sigma} \\
    &= -\int_{\Gamma_f(t)}(\h{\Tt}_f(\h{v}_f\cdot\h{n}_t)^+ + \h{f}_{in}(\h{v}_f\cdot\h{n}_t)^-)\h{\f}_f\,d\h{\sigma}.
\end{align*}
The weak formulation \eqref{eq:transport_weak_ref} can now be derived by writing all spatial derivatives with respect to $x$ instead of $\h{x}$ and performing a change of variables in all integrals. For the surface integral on the right hand side, we use that on $\Gamma_f$, it holds $n_t=n$, $F^{-T}n=n$ (see \eqref{eq:F_representation}), and $v_b\cdot n=0$. Moreover, we employ the identity $\Div(Av_{b/s})=\pa_t J$, which follows from Euler's expansion formula, see e.g.~\cite[Satz 5.2]{eck_mathematische_2017}, and the Piola identity \eqref{eq:piola_identity}.

\begin{defi} \label{defi:weak_formulation}
    A \textbf{weak solution to the fully coupled system} is a tupel $(c,w,q,\Tt_f,\Tt_s)$ satisfying
    \begin{equation} \label{eq:solution_spaces}
        \begin{aligned}
            &c\in C_1^2([0,T]\times[0,L]),& &\\
            &w\in C([0,T],H^1_\Sigma(\Omega_f)^3),&\qquad
            &q\in C([0,T],L^2(\Omega_f)),\\
            &\Tt_f\in L^2(0,T;H^1(\Omega_f)),& 
            &\pa_t(J\Tt_f)\in L^2(0,T;H^1(\Omega_f)'),\\
            &\Tt_s\in L^2(0,T;H^1_{\Gamma_D}(\Omega_s)),&
            &\pa_t(J\Tt_s)\in L^2(0,T;H^1_{\Gamma_D}(\Omega_s)'),
        \end{aligned}
    \end{equation}
    and the ODE:
    \begin{equation} \label{eq:ode_ref}
        \begin{aligned}
            \pa_t c(t,x_1) &= -k\,c(t,x_1) + G(x_1,\mathcal{T}(\Tt_s)(t)),\\
            c(0,x_1) &= c^0(x_1),
        \end{aligned}
    \end{equation}
    for all $(t,x_1)\in[0,T]\times[0,L]$, where
    \begin{equation} \label{eq:T_representation}
        \mathcal{T}(\Tt_s) = K_\gamma\star\mathcal{T}_1(\Tt_s)\qquad\text{with}\qquad\mathcal{T}_1(\Tt_s)(t) = \begin{cases}
            \frac{1}{|\Omega_s|}\int_{\Omega_s}\Tt_s^0(x)\dx &\text{if }t\in(-\gamma,0], \\
            \frac{1}{|\Omega_s|}\Os\Tt_s(t,x)\dx &\text{if }t\in(0,T),
        \end{cases}
    \end{equation}
    the Stokes system:
    \begin{equation} \label{eq:stokes_weak_ref}
        \begin{aligned}
            \Of 2\mu J\,e_F(w):e_F(\phi) - q\,\Div(A\phi)\dx &= -\Of(A^T\G f_b)\cdot\phi + 2\mu J\,e_F(v_b):e_F(\phi)\dx,\\
            -\Of\Div(Aw)\,\psi\dx &= \Of\Div(Av_b)\,\psi\dx,
        \end{aligned}
    \end{equation}
    for all $(\phi,\psi)\in H^1_\Sigma(\Omega_f)^3\times L^2(\Omega_f)$ and in $[0,T]$, and the advection-diffusion equation:
    \begin{equation} \label{eq:transport_weak_ref}
        \begin{aligned}
            \langle\pa_t(J\Tt_f),\f_f\rangle_{H^1(\Omega_f)',H^1(\Omega_f)} + \langle\pa_t(J\Tt_s),\f_s\rangle_{H^1_{\Gamma_D}(\Omega_s)',H^1_{\Gamma_D}(\Omega_s)}~+&\\
            + \int_{\Omega_f} (K^F_f\G\Tt_f \,-\,\Tt_f Aw)\cdot\G\f_f\,dx + \int_{\Omega_s} K^F_s\G\Tt_s\cdot\G\f_s\,dx~+&\\
            + \int_{\Gamma_f}\Tt_f(w\cdot n)^+\f_f J\,d\sigma + \int_\Sigma \alpha(\Tt_f-\Tt_s)(\f_f-\f_s)J|F^{-T}n|\,d\sigma &= -\int_{\Gamma_f}f_{in} (w\cdot n)^-\f_f J\,d\sigma,\\
            (\Tt_f(0),\Tt_s(0)) &= (\Tt_f^0,\Tt_s^0),
        \end{aligned}
    \end{equation}
    for all $(\f_f,\f_s)\in H^1(\Omega_f)\times H^1_{\Gamma_D}(\Omega_s)$ and a.e.~in $(0,T)$, where the radius $R$ is given by:
    \begin{equation} \label{eq:R_representation}
        R(t,x_1):=H(c(t,x_1)),
    \end{equation}
    the deformation $S$ by:
    \begin{equation} \label{eq:S_representation}
        S(t,x)=\begin{cases}\begin{pmatrix}
    	       x_1 \\
    	       \rho(R(t,x_1),\EN{\bar{x}})\frac{\bar{x}}{\EN{\bar{x}}}
            \end{pmatrix}&\qquad\text{if }|\bar{x}|\neq 0, \\
            x &\qquad\text{if }|\bar{x}|=0,
        \end{cases}
    \end{equation}
    and the coefficients depending on $S$ by:
    \begin{equation} \label{eq:coefficients}
        \begin{gathered}
            v_b=\pa_t S, \quad F=DS, \quad J=\det(F), \quad A=JF^{-1},\\
            e_F(w)=\frac{1}{2}\left(F^{-T}\G w + (F^{-T}\G w)^T\right), \quad K_f^F=JF^{-1}K_fF^{-T}, \quad K_s^F=JF^{-1}K_sF^{-T},
        \end{gathered}
    \end{equation}
    and $n$ denotes the outer normal to $\pa\Omega_f$ (on $\Sigma$) and to $\pa\Omega$.
\end{defi}

\begin{rema} \label{rema:homogeneous_dirichlet}
    For simplicity, we choose a homogeneous Dirichlet boundary condition for the solution of the advection-diffusion equation in the solid subdomain, i.e., $\Tt_s\in L^2(0,T;H^1_{\Gamma_D}(\Omega_s))$. The nonhomogeneous Dirichlet data can be included by standard arguments.
\end{rema}

\renewcommand{\labelenumi}{(A\arabic{enumi})}
\renewcommand{\labelenumii}{(A\arabic{enumi}.\arabic{enumii})}
\begin{assu}
    We make the following assumptions about the data and coupling functions in the model \eqref{eq:solution_spaces}--\eqref{eq:coefficients}:
    \begin{enumerate}
        \item The function $\rho:[R_1,R_2]\times[0,\infty)\to[0,\infty)$ satisfies:
        \begin{enumerate}
            \item $\rho(R,R_0)=R$ for all $R\in[R_1,R_2]$,
            \item $\rho(R,r)=r$ for all $(R,r)\in[R_1,R_2]\times\big([0,\infty)\setminus(R_1-3\gamma,R_2+3\gamma)\big)$,
            \item $\rho|_{[R_1,R_2]\times[0,\frac{1}{2}]}\in C^3([R_1,R_2]\times[0,\frac{1}{2}])$, and
            \item $\pa_r\rho(R,r)\geq c_\rho>0$ for all $(R,r)\in[R_1,R_2]\times[0,\infty)$.
        \end{enumerate}
        \item We have $G\in C^3([0,L]\times\R)$ and there exists a constant $C_G>0$ such that
        \begin{equation*}
            \sup_{(x_1,y)\in[0,L]\times\R}|\pa_{x_1}^{\beta_1}\pa_y^{\beta_2}G(x_1,y)|\leq C_G\qquad\text{for all  }\beta\in\NN_0^2\text{ with }0\leq |\beta|\leq 3.
        \end{equation*}
        \item The function $H\in C^3(\R)$ satisfies $R_1\leq H(y) \leq R_2$ for all $y\in\R$ and there exists a constant $C_H>0$ such that
        \begin{equation*}
            \sup_{y\in\R}|\pa_y^{\beta}H(y)|\leq C_H\qquad\text{for }\beta = 1,2,3.
        \end{equation*}
        \item The initial data satisfies $(c^0,\Tt_f^0,\Tt_s^0)\in C^2([0,L])\times L^2(\Omega_f)\times L^2(\Omega_s)$.
        \item The extended boundary data for the normal stress at $\Gamma_f$ fulfills $f_b\in C([0,T],H^1(\Omega_f))$ and the boundary data for the heat flux at $\Gamma_f$ fulfills $f_{in}\in L^2(0,T;L^2(\Gamma_f))$.
        \item The domain parameters are the minimal and maximal radius $0<R_1<R_2<\frac{1}{2}$, the reference radius $R_0\in[R_1,R_2]$, some $\delta>0$ with $R_2+3\delta<\frac{1}{2}$ and $R_1-3\delta>0$, and the length $L>0$.
        \item With the length of the time average interval $\gamma>0$, the kernel $K_\gamma\in C_c(\R)$ is a non-negative function with $\supp(K_\gamma)\subset[0,\gamma)$ and $\int_{\R}K_\gamma(t)\,dt = 1$.
        \item As parameters in the equations, we have the degradation rate constant $k>0$ in the ODE, the dynamic viscosity $\mu>0$ of the fluid, the thermal conductivities $K_f,K_s\in\R^{3\times 3}$ being symmetric, positive definite matrices, the heat transfer coefficient $\alpha>0$ on $\Sigma$, and the final time $T>0$.
    \end{enumerate}
    \label{assu:data}
\end{assu}

\begin{exam} \label{exam:rho}
    A function $\rho$ satisfying Assumption \ref{assu:data} (A1) is given by:
    \begin{align*}
        &\rho:[R_1,R_2]\times [0,\infty)\to[0,\infty),\qquad\rho(R,r)=\int_\R \ov{\rho}(R,s)\,\omega\left(\frac{r-s}{\delta}\right)\,ds,\text{ where}\\
        &\omega:\R\to\R,\qquad\omega(\xi)=
        \begin{cases}
            \left(\int_\R \exp\left(\frac{-1}{1-\zeta^2}\right)\,d\zeta\right)^{-1}\exp\left(\frac{-1}{1-\xi^2}\right) &\text{if }|\xi|<1,\\
            0 &\text{if }|\xi|\geq 1,\text{ and}
        \end{cases}\\
        &\ov{\rho}:[R_1,R_2]\times\R\to\R,\\
        &\qquad\ov{\rho}(R,r)=
        \begin{cases}
            r &\text{if }r\leq R_1-2\delta,\\
            \frac{R-R_1+\delta}{R_0-R_1+\delta}(r-(R_1-2\delta))+R_1-2\delta &\text{if }R_1-2\delta\leq r\leq R_0-\delta,\\
            (r-R_0)+R &\text{if }R_0-\delta\leq r\leq R_0+\delta,\\
            \frac{R_2-R+\delta}{R_2-R_0+\delta}(r-(R_2+2\delta))+R_2+2\delta &\text{if }R_0+\delta\leq r\leq R_2+2\delta,\\
            r &\text{if }R_2+2\delta\leq r,
        \end{cases}
    \end{align*}
    see \cite[Section 4.2, pp.~137-139]{wiedemann_analytical_2023} for more details.
    \qed
\end{exam}

%% file: 3_result.tex
%Result
\section{Global-in-time existence and uniqueness} \label{sec:result}

The main theorem of this paper is the following.

\begin{theo}[Global-in-time existence and uniqueness] \label{theo:main}
    Under the Assumption \ref{assu:data}, there exists a unique weak solution $(c,w,q,\Tt_f,\Tt_s)$ to the fully coupled system \eqref{eq:solution_spaces}--\eqref{eq:coefficients}.
\end{theo}

\begin{figure}[tbh]
    \centering
    \includegraphics[width=0.45\linewidth]{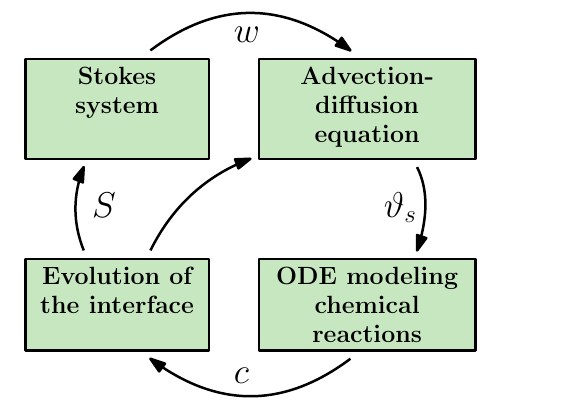} \hfill
    \includegraphics[width=0.45\linewidth]{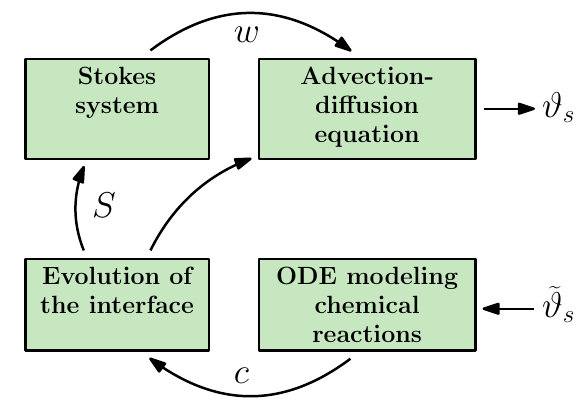}
    \caption{The fully coupled system of equations (left) and the linearized system (right).}
    \label{fig:coupling}
\end{figure}

For the proof of Theorem \ref{theo:main}, we first linearize the system \eqref{eq:solution_spaces}--\eqref{eq:coefficients} and prove in Sections \ref{subsec:c}, \ref{subsec:w_q} and \ref{subsec:theta_f_theta_s} that for any $\widetilde{\Tt}_s\in L^2(0,T;L^2(\Omega_s))$, the subproblems \eqref{eq:ode_ref}, \eqref{eq:stokes_weak_ref} and \eqref{eq:transport_weak_ref} admit a unique solution, see Figure \ref{fig:coupling} for a visualization. More properties of the deformation of the domain can be found in Section \ref{subsec:S}. \par
Then, we can define the solution operator $\mathcal{F}$ of the linearized system:
\begin{equation} \label{eq:solution_operator}
    \mathcal{F}:L^2(0,T;L^2(\Omega_s))\to L^2(0,T;L^2(\Omega_s)),\qquad\mathcal{F}(\widetilde{\Tt}_s)=\Tt_s,
\end{equation}
where $\Tt_s$ is the unique solution of \eqref{eq:transport_weak_ref} given by Proposition \ref{prop:theta_solution_estimates}. In Section \ref{subsec:fixed_point_method}, we use Schaefer's fixed-point theorem to show that $\mathcal{F}$ has a fixed point. Finally, uniqueness is proven in Section \ref{subsec:uniqueness}. \par
When working with the linearized system, a main difficulty is to ensure that the constants only depend on the data from Assumption \ref{assu:data}. In particular, the constants must not depend on the solution of another subproblem. We will pay special attention to this aspect.

%% file: 4_proof.tex
%Proof
\subsection{Existence of a unique solution for the ODE} \label{subsec:c}

\begin{lemm} \label{lemm:T_properties}
    There exists a constant $C>0$, only depending on the data from Assumption \ref{assu:data}, such that for all $\widetilde{\Tt}_s\in L^2(0,T;L^2(\Omega_s))$, the operator $\mathcal{T}$ given by \eqref{eq:T_representation} satisfies $\mathcal{T}(\widetilde{\Tt}_s)\in C([0,T])$ and
    \begin{equation} \label{eq:T_estimates}
        \N{\mathcal{T}(\widetilde{\Tt}_s)}_{C([0,T])} \leq C\,\left(\N{\Tt_s^0}_{\Lt(\Omega_s)} + \N{\widetilde{\Tt}_s}_{\Lt(0,T;\Lt(\Omega_s))}\right).
    \end{equation}
\end{lemm}

\begin{proo}
    The first claim follows from the convolution with $K_\gamma\in C_c(\R)$. For the estimate, we have
    \begin{align*}
        \N{\mathcal{T}_1(\widetilde{\Tt}_s)}_{L^1(-\gamma,T)} &= \frac{1}{|\Omega_s|}\int_{-\gamma}^0\EN{\int_{\Omega_s}\Tt_s^0(x)\dx}\dt + \frac{1}{|\Omega_s|}\int_0^T\EN{\Os\widetilde{\Tt}_s(t,x)\dx}\dt \\
        &\leq \frac{\gamma}{\sqrt{|\Omega_s|}}\N{\Tt_s^0}_{L^2(\Omega_s)} + \frac{\sqrt{T}}{\sqrt{|\Omega_s|}}\N{\widetilde{\Tt}_s}_{L^2(0,T;L^2(\Omega_s))},
    \end{align*}
    and thus \eqref{eq:T_estimates} by Young's convolution inequality, where $C$ depends in particular on the choice of $K_\gamma$ in Assumption \ref{assu:data} (A7).
    \qed
\end{proo}

\begin{prop} \label{prop:c_solution_estimates}
    There exists a constant $C>0$, only depending on the data from Assumption \ref{assu:data}, such that for all $\widetilde{\Tt}_s\in L^2(0,T;L^2(\Omega_s))$, there exists a unique solution $c\in C_1^2([0,T]\times[0,L])$ of the ODE \eqref{eq:ode_ref} and \eqref{eq:T_representation}, given by
    \begin{equation} \label{eq:c_solution}
        c(t,x_1) = e^{-kt}\,c^0(x_1) + e^{-kt}\int_0^t e^{ks}\,G(x_1,\mathcal{T}(\widetilde{\Tt}_s)(s))\,ds,
    \end{equation}
    and satisfying
    \begin{equation} \label{eq:c_estimates}
        \sup_{(t,x_1)\in[0,T]\times[0,L]} |\pa_t^{\beta_1}\pa_{x_1}^{\beta_2}c(t,x_1)| \leq C\qquad\text{for all }\beta\in\NN_0^2\text{ with }0\leq |\beta|\leq 2\text{ and }\beta_1\leq 1.
    \end{equation}
\end{prop}

\begin{proo}
    From Lemma \ref{lemm:T_properties} and Assumption \ref{assu:data} (A2), it follows that $(t,x_1)\mapsto G(x_1,\mathcal{T}(\widetilde{\Tt}_s)(t))$ and its first and second derivative with respect to $x_1$ are in $C([0,T]\times[0,L])$, and thus we obtain existence, uniqueness and \eqref{eq:c_solution} from the variation of constants formula, see e.g.~\cite[p.~19]{wilke_gewohnliche_2011}. The estimates \eqref{eq:c_estimates} result from \eqref{eq:c_solution} and the explicit representation
    \begin{equation} \label{eq:c_time_derivative}
        \pa_t c(t,x_1) = -k\,e^{-kt}\,c^0(x_1) - k\,e^{-kt}\int_0^t e^{ks}\,G(x_1,\mathcal{T}(\widetilde{\Tt}_s)(s))\,ds + G(x_1,\mathcal{T}(\widetilde{\Tt}_s)(t)),
    \end{equation}
    and from the derivatives of \eqref{eq:c_solution} and \eqref{eq:c_time_derivative} with respect to $x_1$, using Assumption \ref{assu:data} (A2) and (A4).
    \qed
\end{proo}

From Proposition \ref{prop:c_solution_estimates} and Assumption \ref{assu:data} (A3), we deduce the

\begin{coro} \label{coro:R_solution_estimates}
    There exists a constant $C>0$, only depending on the data from Assumption \ref{assu:data}, such that for all $\widetilde{\Tt}_s\in L^2(0,T;L^2(\Omega_s))$, the radius $R$ given by \eqref{eq:ode_ref}, \eqref{eq:T_representation}, and \eqref{eq:R_representation} satisfies $R\in C_1^2([0,T]\times[0,L])$, $R_1\leq R(t,x_1)\leq R_2$ for all $(t,x_1)\in[0,T]\times[0,L]$ and
    \begin{equation} \label{eq:R_estimates}
        \sup_{(t,x_1)\in[0,T]\times[0,L]} |\pa_t^{\beta_1}\pa_{x_1}^{\beta_2}R(t,x_1)| \leq C\qquad\text{for all }\beta\in\NN_0^2\text{ with }1\leq |\beta|\leq 2\text{ and }\beta_1\leq 1.
    \end{equation}
\end{coro}

\subsection{Deformation of the domain and its properties} \label{subsec:S}

For convenience, we introduce the subdomain 
\begin{equation} \label{eq:Z_definition}
    \textstyle Z:=\left\{(x_1,\bar{x})\in\Omega: R_1-3\delta<|\bar{x}|<R_2+3\delta\right\}.
\end{equation}

\renewcommand{\labelenumi}{(\arabic{enumi})}
\begin{rema} \label{rema:S_bijectivity}
    Let us summarize some direct consequences of Assumption \ref{assu:data}. Choose any $\widetilde{\Tt}_s\in L^2(0,T;L^2(\Omega_s))$, let $R\in C_1^2([0,T]\times[0,L])$ be given by Corollary \ref{coro:R_solution_estimates} and $S$ by \eqref{eq:S_representation}.
    \begin{enumerate}
        \item From (A1.2), (A1.3) and (A1.4), it follows that for any $R^\ast\in[R_1,R_2]$, the function $\rho(R^\ast,\cdot):[0,\infty)\to[0,\infty)$ is strictly increasing.
        \item By (A1.2) and (A1.3), we have $\sup_{(R^\ast,r)\in[R_1,R_2]\times[0,\infty)}|\pa^{\beta}\rho(R^\ast,r)|<\infty$ for all $\beta\in\NN_0^2$ with $1\leq|\beta|\leq 3$.
        \item Assumption (A1.2) implies $S(t,x)=x$ for all $t\in[0,T]$ and all $x\in\ov{\Omega\setminus Z}$. Moreover, by its definition \eqref{eq:S_representation}, $S(t,\cdot)$ does not deform $\ov{\Omega}$ in the $x_1$-direction and only radially in each $x_2$-$x_3$-section. As $\rho(R^\ast,\cdot)$ is strictly increasing for any $R^\ast\in[R_1,R_2]$, the deformation $S(t,\cdot):\ov{\Omega}\to\ov{\Omega}$ is one-to-one.
        \item For the sake of clarity, we provide the following representations of the evolving configuration. By (A1.1),
        \begin{align*}
            \Sigma(t) &= S(t,\Sigma) = \left\{(x_1,\bar{x})\in\ov{\Omega}: |\bar{x}|=R(t,x_1)\right\},\\ 
            \intertext{and in this sense, $R$ is the radius of the evolving interface. As $\rho(R^\ast,\cdot)$ is strictly increasing for any $R^\ast\in[R_1,R_2]$, we obtain:}
            \Omega_f(t) &= S(t,\Omega_f) = \left\{(x_1,\bar{x})\in\Omega: |\bar{x}|<R(t,x_1)\right\},\text{ and}\\
            \Omega_s(t) &= S(t,\Omega_s) = \Omega\setminus\overline{\Omega_f(t)}.\\
            \intertext{For the boundaries, it holds:}
            \Gamma_f(t) &= S(t,\Gamma_f) = \left\{(x_1,\bar{x})\in\ov{\Omega}: x_1\in\{0,L\}\text{ and }|\bar{x}|\leq R(t,x_1)\right\},\text{ and}\\
            \Gamma_D(t) &= S(t,\Gamma_D) = \Gamma_D.
        \end{align*}
    \end{enumerate}
\end{rema}

In the following, we use the notation \eqref{eq:coefficients}, in particular $v_b=\pa_t S$, $F=DS$, $J=\det(F)$, $A=JF^{-1}$, $K_{f/s}^F=JF^{-1}K_{f/s}F^{-T}$, where `$f/s$' indicates an equation for both `$f$' and `$s$'.
    
\begin{prop} \label{prop:S_solution_estimates}
    There exist constants $C,\eta, c_f^K, c_s^K >0$, only depending on the data from Assumption \ref{assu:data}, such that for all $\widetilde{\Tt}_s\in L^2(0,T;L^2(\Omega_s))$, the deformation $S$ given by \eqref{eq:ode_ref}, \eqref{eq:T_representation}, \eqref{eq:R_representation}, and \eqref{eq:S_representation} satisfies $S\in C_1^2([0,T]\times\ov{\Omega})^3$ and 
    \begin{align}
        &\sup_{(t,x)\in[0,T]\times\ov{\Omega}}|\pa_t^{\beta_1}\pa_{x_1}^{\beta_2}\pa_{x_2}^{\beta_3}\pa_{x_3}^{\beta_4}S(t,x)|\leq C\qquad\text{for all  }\beta\in\NN_0^4\text{ with }0\leq |\beta|\leq 2\text{ and }\beta_1\leq 1,\label{eq:S_estimate_derivatives}\\
        &\inf_{(t,x)\in[0,T]\times\ov{\Omega}}J(t,x) \geq \eta,\label{eq:S_estimate_J}\\
        &\inf_{(t,x)\in[0,T]\times\Sigma}J(t,x)\,|F^{-T}(t,x)n(x)| \geq \frac{R_1}{R_0},\label{eq:S_estimate_interface}\\
        &\inf_{(t,x)\in[0,T]\times\ov{\Omega}}\xi^T K_{f/s}^F(t,x)\xi \geq c^K_{f/s}\,|\xi|^2\qquad\text{for all }\xi\in\R^3,\label{eq:S_estimate_K}\\
        &\Div(A)=0\qquad\text{in }[0,T]\times\ov{\Omega},\qquad\text{(Piola identity)}\label{eq:piola_identity}
    \end{align}
    where $n$ denotes the outer normal to $\Omega_f$. In particular, we have:
    \begin{align}
        \N{F}_{C^1([0,T]\times\ov{\Omega})^{3\times 3}} + \N{J}_{C^1([0,T]\times\ov{\Omega})} + \N{A}_{C^1([0,T]\times\ov{\Omega})^{3\times 3}} &\leq \widetilde{C}(C,\eta), \label{eq:coefficients_estimates} \\
        \N{v_b}_{C([0,T]\times\ov{\Omega})^3} + \N{\G v_b}_{C([0,T]\times\ov{\Omega})^{3\times 3}} &\leq C. \label{eq:v_b_estimates}
    \end{align}
\end{prop}

\begin{proo}
    First, remark that $S\in C_1^2([0,T]\times\ov{\Omega})^3$. Indeed, on the one hand, Assumption \ref{assu:data} (A1.2) implies 
    \begin{equation} \label{eq:S_identity_map}
        S(t,x) = x\qquad\text{for all }(t,x)\in[0,T]\times\ov{\Omega\setminus Z},
    \end{equation}
    and on the other hand,
    \begin{equation*}
        S\in C_1^2([0,T]\times\ov{Z})^3,
    \end{equation*}
    because in $\ov{Z}$, $S$ is the composition of sufficiently smooth functions, see \eqref{eq:S_representation}, Assumption \ref{assu:data} (A1.3), and Corollary \ref{coro:R_solution_estimates}.\par
    Due to \eqref{eq:S_identity_map}, we have $F\equiv\Id_3$, $J\equiv 1$, and $\pa_t S\equiv 0$ in $[0,T]\times\ov{\Omega\setminus Z}$, so it suffices to prove \eqref{eq:S_estimate_derivatives} and \eqref{eq:S_estimate_J} for $\ov{Z}$ instead of $\ov{\Omega}$.\par
    To derive the estimate \eqref{eq:S_estimate_derivatives}, first compute the derivatives of $S$ in $\ov{Z}$ from $\eqref{eq:S_representation}$ by the chain rule. Then, use Assumption \ref{assu:data} (A1) and Corollary \ref{coro:R_solution_estimates} to obtain \eqref{eq:S_estimate_derivatives} in $[0,T]\times\ov{Z}$.\par
    For \eqref{eq:S_estimate_J} and \eqref{eq:S_estimate_interface}, note that in $[0,T]\times\ov{Z}$, $F$ and $F^{-T}$ are of the forms
    \begin{equation} \label{eq:F_representation}
        F(t,x) = \begin{pmatrix}1 & 0 \\ \mathcal{A}(t,x) & \mathcal{B}(t,x) 
        \end{pmatrix}\qquad\text{and}\qquad F^{-T}(t,x) = \begin{pmatrix}1 & -\mathcal{A}^T(t,x)\mathcal{B}^{-T}(t,x) \\ 0 & \mathcal{B}^{-1}(t,x)
        \end{pmatrix}
    \end{equation}
    with $\mathcal{A}(t,x)\in\R^{2\times 1}$ and $\mathcal{B}(t,x)=\mathcal{B}^T(t,x)\in\R^{2\times 2}$, so after some computation, we obtain
    \begin{equation*}
        J(t,x) = \det(\mathcal{B})(t,x_1,\bar{x}) = \frac{1}{|\bar{x}|}\,\rho(R(t,x_1),|\bar{x}|)\,\pa_r\rho(R(t,x_1),|\bar{x}|) \geq \frac{R_1-3\delta}{R_2+3\delta}\,c_\rho =: \eta > 0
    \end{equation*}
    for all $(t,x)\in[0,T]\times\ov{Z}$ -- which justifies taking the inverse of $F(t,x)$ -- and
    \begin{equation*}
        J(t,x)\,|F^{-T}(t,x)n(x)| = \frac{R(t,x_1)}{R_0}\sqrt{(\pa_{x_1}R(t,x_1))^2 + 1} \geq \frac{R_1}{R_0}
    \end{equation*}
    for all $(t,x)\in[0,T]\times\Sigma$. \par
    To prove \eqref{eq:S_estimate_K}, fix any $\xi\in\R^3$. By Assumption \ref{assu:data} (A8), $K_{f/s}$ are symmetric, positive definite matrices, so denoting by $\lambda_{f/s}$ the smallest eigenvalue of $K_{f/s}$, we get, in $[0,T]\times\ov{\Omega}$:
    \begin{equation*}
        \xi^T K_{f/s}^F\xi = \xi^T(JF^{-1}K_{f/s}F^{-T})\xi = J(F^{-T}\xi)^TK_{f/s}(F^{-T}\xi) \geq \eta\,\lambda_{f/s}\,\xi^TF^{-1}F^{-T}\xi \geq 0.
    \end{equation*}
    As $F^{-1}F^{-T}$ is also symmetric and positive definite, we can denote its smallest eigenvalue by $\lambda_1(t,x)>0$, and by $\lambda_2(t,x)>0$ and $\lambda_3(t,x)>0$ the other two, such that in $[0,T]\times\ov{\Omega}$:
    \begin{align*}
        &\lambda_1\lambda_2\lambda_3 = \det(F^{-1}F^{-T}) = \frac{1}{J^2} \geq \frac{1}{\widetilde{C}},\\
        &\lambda_1+\lambda_2+\lambda_3 = \trc(F^{-1}F^{-T}) = F^{-T}:F^{-T} \leq \widetilde{C},
    \end{align*}
    where $\widetilde{C}>0$ is a polynomial in the constants $C$ and $\frac{1}{\eta}$ from \eqref{eq:S_estimate_derivatives} and \eqref{eq:S_estimate_J}. Combining these estimates, we end up with
    \begin{equation*}
        \xi^T K_{f/s}^F\xi \geq \eta\,\lambda_{f/s}\,\xi^TF^{-1}F^{-T}\xi \geq \eta\,\lambda_{f/s}\,\lambda_1\,|\xi|^2 \geq \eta\,\lambda_{f/s}\,\frac{1}{\widetilde{C}\lambda_2\lambda_3}\,|\xi|^2 \geq \eta\,\lambda_{f/s}\,\frac{1}{\widetilde{C}^3}\,|\xi|^2.
    \end{equation*}
    in $[0,T]\times\ov{\Omega}$, so we can define $c_{f/s}^K := \eta\,\lambda_{f/s}\,\frac{1}{\widetilde{C}^3}$. \par
    For the proof of the Piola identity, see e.g.~\cite[p.~39]{ciarlet_mathematical_1988}.
    \qed
\end{proo}

\begin{coro} \label{coro:Lipschitz_domain}
    For all $\widetilde{\Tt}_s\in L^2(0,T;L^2(\Omega_s))$ and all $t\in[0,T]$, the domains $\Omega_f(t)=S(t,\Omega_f)$ and $\Omega_s(t)=S(t,\Omega_s)$ are Lipschitz domains, where $S$ is given by \eqref{eq:ode_ref}, \eqref{eq:T_representation}, \eqref{eq:R_representation}, and \eqref{eq:S_representation}.
\end{coro}

\subsection{Existence of a unique solution for the Stokes system} \label{subsec:w_q}

In this section, we write $\V$ for $H_\Sigma^1(\Omega_f)^3$ and $\Q$ for $L^2(\Omega_f)$, because the Stokes system is set only in $\Omega_f$. We denote the dual pairings by $\DP{\cdot,\cdot}{\V}$ and $\DP{\cdot,\cdot}{\Q}$ and the $L^2$-inner product by $(\cdot,\cdot)_\Q$. As norms, choose $\N{q}_\Q = \N{q}_{L^2(\Omega_f)}$ and $\N{v}_\V = \N{\G v}_{L^2(\Omega_f)}$, which is a norm on $\V$ by Poincaré's inequality. We write $\Bil(\V,\Q)$ for the space of continuous bilinear forms on $\V\times \Q$ with the usual operator norm
\begin{equation*}
    \N{f}_\Bil := \sup_{v\in \V\setminus\{0\},~q\in \Q\setminus\{0\}}\frac{|f(v,q)|}{\N{v}_\V\N{q}_\Q},
\end{equation*}
and we are particularly interested in the following subsets:
\begin{align*}
    &A_C^{\alpha_0} := \{a\in\Bil(\V,\V): \N{a}_\Bil\leq C\text{ and }a(v,v)\geq\alpha_0\N{v}_\V^2\text{ for all }v\in \V\},\\
    &B^{\beta_0} := \left\{b\in\Bil(\V,\Q): \inf_{q\in \Q\setminus\{0\}}\sup_{v\in \V\setminus\{0\}}\frac{b(v,q)}{\N{v}_\V\N{q}_\Q}\geq\beta_0\right\},\\
    &V'_C := \{f\in \V': \N{f}_{\V'}\leq C\},\\
    &Q'_C := \{g\in \Q': \N{g}_{\Q'}\leq C\}.
\end{align*}
Now, for all $t\in[0,T]$, define the (bi)linear forms
\begin{equation} \label{eq:stokes_bilinear_forms}
    \begin{aligned}
        &a(t): \V\times \V\to\R,\qquad a(t)(w,\phi) = \Of 2\mu J(t)\,e_{F(t)}(w):e_{F(t)}(\phi)\dx,\\
        &b(t): \V\times \Q\to\R,\qquad b(t)(\phi,q) = -\Of q\,\Div(A(t)\phi)\dx = -\Of q\,A(t):\G\phi\dx, \\
        &f(t): \V\to\R,\qquad \DP{f(t),\phi}{\V} = -\Of(A^T(t)\G f_b(t))\cdot\phi\dx - a(t)(v_b(t),\phi),\\
        &g(t): \Q\to\R,\qquad \DP{g(t),\psi}{\Q} = \Of \Div(A(t)v_b(t))\,\psi\dx = -b(t)(v_b(t),\psi),
    \end{aligned}
\end{equation}
where the reformulation of $b(t)$ follows from the product rule and the Piola identity \eqref{eq:piola_identity}. Using this notation, the Stokes system \eqref{eq:stokes_weak_ref} can be written as
\begin{equation} \label{eq:stokes_weak_general_form}
    \left\{\begin{aligned}
        &a(t)(w(t),\phi) + b(t)(\phi,q(t)) = \DP{f(t),\phi}{\V},\\
        &b(t)(w(t),\psi) = \DP{g(t),\psi}{\Q},
    \end{aligned}\right.
\end{equation}
for all $(\phi,\psi)\in \V\times \Q$ and all $t\in[0,T]$. The proof of the following proposition relies on two general results for problems of the form \eqref{eq:stokes_weak_general_form}, which can be found, for example, in \cite[Corollary 4.2.1]{boffi_mixed_2013} and \cite[Lemma B.4]{gahn_rigorous_2025}. This method has been applied before, e.g., in \cite{wiedemann_homogenisation_2024}.
\begin{prop} \label{prop:w_solution_estimates}
    There exist constants $\alpha_0,\beta_0,C>0$, only depending on the data from Assumption \ref{assu:data}, such that for all $\widetilde{\Tt}_s\in L^2(0,T;L^2(\Omega_s))$, the (bi)linear forms \eqref{eq:stokes_bilinear_forms} with coefficients given by \eqref{eq:ode_ref}, \eqref{eq:T_representation}, \eqref{eq:R_representation}, \eqref{eq:S_representation}, and \eqref{eq:coefficients} satisfy
    \begin{equation} \label{eq:abfg_solution}
        (a(t),b(t),f(t),g(t))\in A_C^{\alpha_0}\times B^{\beta_0}\times V'_C\times Q'_C\qquad\text{for all }t\in[0,T],
    \end{equation}
    and the Stokes system \eqref{eq:stokes_weak_ref} admits a unique solution $(w,q)\in C([0,T],\V)\times C([0,T],\Q)$ with
    \begin{equation} \label{eq:w_estimates}
        \sup_{t\in[0,T]}\N{w(t)}_{\V}\leq\frac{C}{\alpha_0}+\frac{2C^\frac{3}{2}}{\alpha_0^\frac{1}{2}\beta_0}\qquad\text{and}\qquad\sup_{t\in[0,T]}\N{q(t)}_{\Q}\leq\frac{2C^\frac{3}{2}}{\alpha_0^\frac{1}{2}\beta_0}+\frac{C^2}{\beta_0^2}.
    \end{equation}
\end{prop}

\begin{proo}
    Choose any $\widetilde{\Tt}_s\in L^2(0,T;L^2(\Omega_s))$ and let the coefficients for \eqref{eq:stokes_weak_ref} and \eqref{eq:stokes_bilinear_forms} be given by \eqref{eq:ode_ref}, \eqref{eq:T_representation}, \eqref{eq:R_representation}, \eqref{eq:S_representation}, and \eqref{eq:coefficients}. First, fix any time $t\in[0,T]$. Then, for any $w,\phi\in \V$ and any $q,\psi\in \Q$,
    \begin{align*}
        &\begin{aligned}
            |a(t)(w,\phi)| &\leq \frac{\mu}{2}\N{J(t)}_{C(\ov{\Omega_f})}\N{F^{-T}(t)\G w + (F^{-T}(t)\G w)^{T}}_\Q\N{F^{-T}(t)\G\phi + (F^{-T}(t)\G\phi)^{T}}_\Q \\
            &\leq C \N{w}_\V\N{\phi}_\V,
        \end{aligned}\\
        &|b(t)(\phi,q)| \leq \N{A(t)}_{C(\ov{\Omega_f})}\N{q}_\Q\N{\G\phi}_\Q \leq C \N{q}_\Q\N{\phi}_\V, \\
        &|\DP{f(t),\phi}{\V}| \leq \N{A(t)}_{C(\ov{\Omega_f})}\N{\G f_b(t)}_\Q\N{\phi}_\Q + \widetilde{C}\N{v_b(t)}_\V\N{\phi}_\V \leq C\N{\phi}_\V,\\
        &|\DP{g(t),\psi}{\Q}| \leq \N{A(t)}_{C(\ov{\Omega_f})}\N{v_b(t)}_\V\N{\psi}_\Q \leq C\N{\psi}_\Q,
    \end{align*}
    with some constant $C>0$ which -- due to Proposition \ref{prop:S_solution_estimates} -- only depends on the data from Assumption \ref{assu:data}. Note that we have used Poincaré's inequality $\N{v}_\Q\leq C\N{v}_\V$, and we need \eqref{eq:S_estimate_derivatives} \textsl{and} \eqref{eq:S_estimate_J} to estimate $F^{-T}$. \par
    To show that $a(t)$ is coercive, we use a generalized Korn inequality that can be found, e.g., in \cite[Theorem 3.3]{mielke_thermoviscoelasticity_2020}. By Proposition \ref{prop:S_solution_estimates}, we have constants $C_K,c_K>0$, only depending on the data from Assumption \ref{assu:data}, such that
    \begin{equation*}
        F^{-T}(t)\in\{M\in C^{0,1}(\ov{\Omega_f})^{3\times 3}: \N{M}_{C^{0,1}(\ov{\Omega_f})}\leq C_K\text{ and }\det(M)\geq c_K\}.
    \end{equation*}
    Thus, by the generalized Korn inequality, we get a constant $\alpha_0>0$, only depending on the data from Assumption \ref{assu:data}, such that for all $v\in \V$:
    \begin{equation*}
        a(t)(v,v)\geq \frac{1}{2}\mu\eta\N{F^{-T}(t)\G v + (F^{-T}(t)\G v)^{T}}_\Q^2 \geq \alpha_0\N{v}_\V^2,
    \end{equation*}
    where $\eta$ is from \eqref{eq:S_estimate_J}. \par
    Next, for the inf-sup estimate of $b(t)$, choose any $q\in L^2(\Omega_f)\setminus\{0\}$. By Bogovski\u{i}'s theorem (see e.g. \cite[Theorem 5.4]{fabricius_stokes_2019}), there exists a constant $C_B>0$, only depending on $\Omega_f$ and $\Sigma$, such that we can choose $\tilde{v}\in \V$ with
    \begin{equation*}
        \Div(\tilde{v})=q\text{ in }\Omega_f\qquad\text{and}\qquad\N{\tilde{v}}_{\V}\leq C_B\N{q}_\Q.
    \end{equation*}
    As $J^{-1}(t)F(t)\in C^1(\ov{\Omega_f})$ by Proposition \ref{prop:S_solution_estimates}, $v:=J^{-1}(t)F(t)\tilde{v}$ defines a function $v\in H^1_{\Sigma}(\Omega_f)\setminus\{0\}$, and by the product rule and Poincaré's inequality,
    \begin{equation*}
        \N{v}_\V = \N{\G(J^{-1}(t)F(t)\tilde{v})}_\Q \leq \frac{1}{\beta_0 C_B}\N{\tilde{v}}_\V\leq\frac{1}{\beta_0}\N{q}_\Q
    \end{equation*}
    with some $\beta_0>0$ which -- due to Proposition \ref{prop:S_solution_estimates} -- only depends on the data from Assumption \ref{assu:data}. Consequently, we obtain
    \begin{align*}
        \sup_{w\in \V\setminus\{0\}}\frac{b(t)(w,q)}{\N{w}_\V} &= \sup_{w\in \V\setminus\{0\}}\frac{-(\Div(A(t)w),q)_\Q}{\N{w}_\V} \\
        &\geq \frac{(\Div(A(t)v),q)_\Q}{\N{v}_\V} = \frac{(\Div(\tilde{v}),q)_\Q}{\N{v}_\V} = \frac{\N{q}_\Q^2}{\N{v}_\V} \geq \beta_0\frac{\N{q}_\Q^2}{\N{q}_\Q} = \beta_0\N{q}_\Q.
    \end{align*}\par
    Summing up, we have proven \eqref{eq:abfg_solution}. By a general result for problems of the form \eqref{eq:stokes_weak_general_form} that can be found, for example, in \cite[Corollary 4.2.1]{boffi_mixed_2013}, \eqref{eq:abfg_solution} implies that the Stokes system \eqref{eq:stokes_weak_ref} has a unique weak solution $(w(t),q(t))\in \V\times \Q$ satisfying
    \begin{equation*}
        \N{w(t)}_{\V}\leq\frac{C}{\alpha_0}+\frac{2C^\frac{3}{2}}{\alpha_0^\frac{1}{2}\beta_0}\qquad\text{and}\qquad\N{q(t)}_{\Q}\leq\frac{2C^\frac{3}{2}}{\alpha_0^\frac{1}{2}\beta_0}+\frac{C^2}{\beta_0^2}
    \end{equation*}
    for each $t\in[0,T]$, so in particular, we obtain \eqref{eq:w_estimates}. \par
    To prove continuity in time, we first show that the (bi)linear forms \eqref{eq:stokes_bilinear_forms} are continuous in time. Therefore, choose any $t_0\in[0,T]$ and let $(t_k)_k\subset [0,T]$ be any sequence converging to $t_0$. Then, for all $u,v\in \V$,
    \begin{align*}
        &|a(t_k)(u,v)-a(t_0)(u,v)| \\
        &= \mu\Big|\Of\left((F^{-T}(t_k)-F^{-T}(t_0))\G u + ((F^{-T}(t_k)-F^{-T}(t_0))\G u)^T\right):\left(A^T(t_k)\G v\right)\\
        &\hspace{1.8cm} + \left(F^{-T}(t_0)\G u + (F^{-T}(t_0)\G u)^T\right):\left((A^T(t_k)-A^T(t_0))\G v\right)\dx\Big|\\
        &\leq \mu\widetilde{C}\N{u}_\V\N{v}_\V\left(\N{F^{-T}(t_k)-F^{-T}(t_0)}_{C(\ov{\Omega_f})}\N{A^T(t_k)}_{C(\ov{\Omega_f})}\right.\\
        &\hspace{3.7cm} + \left.\N{F^{-T}(t_0)}_{C(\ov{\Omega_f})}\N{A^T(t_k)-A^T(t_0)}_{C(\ov{\Omega_f})}\right),
    \end{align*}
    where in the first step, we use $J\,e_{F}(u):e_{F}(v) = e_{F}(u):(A^{-T}\G v)$, and $\widetilde{C}$ only depends on the number of components and terms. As $A$ and $F^{-1}$ are continuous in the compact set $[0,T]\times\ov{\Omega_f}$, see Proposition \ref{prop:S_solution_estimates}, we obtain
    \begin{equation*}
        \N{a(t_k)-a(t_0)}_\Bil\to 0\qquad\text{as }k\to\infty.
    \end{equation*}
    Similarly, for all $v\in \V$ and all $p\in \Q$,
    \begin{align*}
        |-b(t_k)(v,p)+b(t_0)(v,p)| &= \Big|\Of p\,(A(t_k)-A(t_0)):\G v\dx\Big| \\
        &\leq \widetilde{C}\N{v}_\V\N{p}_\Q\N{A(t_k)-A(t_0)}_{C(\ov{\Omega_f})},
    \end{align*}
    and thus
    \begin{equation*}
        \N{b(t_k)-b(t_0)}_\Bil\to 0\qquad\text{as }k\to\infty.
    \end{equation*}
    Furthermore, using Poincaré's inequality, for all $\phi\in \V$, we get:
    \begin{align*}
        &|-\DP{f(t_k),\phi}{\V} + \DP{f(t_0),\phi}{\V}|\\
        &=\Big|\Of \left(A^T(t_k)(\G f_b(t_k) - \G f_b(t_0)) + (A^T(t_k)-A^T(t_0))\G f_b(t_0)\right)\cdot\phi\dx \\
        &\hspace{1.8cm}+~a(t_k)(v_b(t_k),\phi) - a(t_0)(v_b(t_k),\phi) + a(t_0)(v_b(t_k)-v_b(t_0),\phi)\Big|\\
        &\leq \N{\phi}_\V\left(\N{A^T(t_k)}_{C(\ov{\Omega_f})}\N{f_b(t_k)-f_b(t_0)}_{H^1} + \N{A^T(t_k)-A^T(t_0)}_{C(\ov{\Omega_f})}\N{f_b(t_0)}_{H^1}\right.\\
        &\hspace{1.8cm} + \left.\mu\N{a(t_k)-a(t_0)}_{\Bil}\N{\G v_b(t_k)}_{C(\ov{\Omega_f})} + \mu\N{a(t_0)}_\Bil\N{\G v_b(t_k) - \G v_b(t_0)}_{C(\ov{\Omega_f})}\right).
    \end{align*}
    By the continuity of $A$ and $\G v_b = \G(\pa_t S)$ in $[0,T]\times\ov{\Omega_f}$ (see Proposition \ref{prop:S_solution_estimates}) and Assumption \ref{assu:data} (A5), it follows
    \begin{equation*}
        \N{f(t_k)-f(t_0)}_{\V'}\to 0\qquad\text{as }k\to\infty.
    \end{equation*}
    Finally, for all $\psi\in \Q$, we can estimate
    \begin{align*}
        |\DP{g(t_k),\psi}{\Q} - \DP{g(t_0),\psi}{\Q}|&=\Big|\Of \psi\,(A(t_k):\G v_b(t_k) - A(t_0):\G v_b(t_0))\dx\Big| \\
        &\leq\widetilde{C}\N{\psi}_\Q\N{A(t_k):\G v_b(t_k) - A(t_0):\G v_b(t_0)}_{C(\ov{\Omega_f})},
    \end{align*}
    and because $A:\G v_b = A:\G(\pa_t S)$ is continuous in $[0,T]\times\ov{\Omega_f}$, it follows
    \begin{equation*}
        \N{g(t_k)-g(t_0)}_{\Q'}\to 0\qquad\text{as }k\to\infty.
    \end{equation*}
    In summary, we have $a\in C([0,T],\Bil(\V,\V))$, $b\in C([0,T],\Bil(\V,\Q))$, $f\in C([0,T],\V')$ and $g\in C([0,T],\Q')$. By the general result of \cite[Lemma B.4]{gahn_rigorous_2025} for problems of the form \eqref{eq:stokes_weak_general_form}, the operator that assigns the solution $(w,q)$ of the Stokes system \eqref{eq:stokes_weak_ref} to the (bi)linear forms $(a,b,f,g)$ is Lipschitz continuous in $A_C^{\alpha_0}\times B^{\beta_0}\times V'_C\times Q'_C$. So, as composition of continuous maps, we obtain $w\in C([0,T],\V)$ and $q\in C([0,T],\Q)$.
    \qed
\end{proo}

\subsection{Existence of a unique solution for the advection-diffusion equation} \label{subsec:theta_f_theta_s}

The following product rule in $L^2(0,T;H^1(\Omega_f)')$ resp. $L^2(0,T;H^1_{\Gamma_D}(\Omega_s)')$ is essential for solving the advection-diffusion equation \eqref{eq:transport_weak_ref}. Again, `$f/s$' indicates an equation for both `$f$' and `$s$'.

\begin{lemm} \label{lemm:theta_product_rule}
    Suppose that for $\Tt_f\in L^2(0,T;H^1(\Omega_f))$ resp. $\Tt_s\in L^2(0,T;H^1_{\Gamma_D}(\Omega_s))$ and $J\in C^1([0,T]\times\ov{\Omega})$ with $J\geq\eta>0$, the generalized time derivative $\pa_t(J\Tt_f)\in L^2(0,T;H^1(\Omega_f)')$ resp. $\pa_t(J\Tt_s)\in L^2(0,T;H^1_{\Gamma_D}(\Omega_s)')$ exists. Then also the generalized time derivatives
    \begin{equation*}
        \pa_t\Tt_f,\pa_t(\sqrt{J}\Tt_f)\in L^2(0,T;H^1(\Omega_f)'),\qquad \pa_t\Tt_s,\pa_t(\sqrt{J}\Tt_s)\in L^2(0,T;H^1_{\Gamma_D}(\Omega_s)')
    \end{equation*}
    exist, and we have the product rules:
    \begin{align}
        \pa_t(J\Tt_{f/s}) &= J\,\pa_t\Tt_{f/s}+\pa_t J\,\Tt_{f/s},\label{eq:product_rule_J}\\
        \pa_t(J\Tt_{f/s}) &= \sqrt{J}\,\pa_t(\sqrt{J}\Tt_{f/s})+\pa_t(\sqrt{J})\sqrt{J}\Tt_{f/s} = \sqrt{J}\,\pa_t(\sqrt{J}\Tt_{f/s})+\frac{1}{2}\pa_t J\,\Tt_{f/s}.\label{eq:product_rule_sqrt_J}
    \end{align}
\end{lemm}

\begin{proo}
    We prove \eqref{eq:product_rule_J} for $\Tt_f\in L^2(0,T;H^1(\Omega_f))$ and $\pa_t(J\Tt_f)\in L^2(0,T;H^1(\Omega_f)')$, the rest follows in the same way. For any $\f_f\in H^1(\Omega_f)$ and any $\f\in C_c^\infty(0,T)$, the generalized integration by parts formula \cite[Proposition 23.23 (iv)]{zeidler_nonlinear_1990} and the classical product rule yield
    \begin{align*}
        &-\int_0^T\DP{\pa_t(J\Tt_f),\frac{1}{J}\f_f}{H^1(\Omega_f)}\f\,dt + \int_0^T\DP{\pa_t J\,\Tt_f,\frac{1}{J}\f_f}{H^1(\Omega_f)}\f\,dt \\
        &=\int_0^T\left(\pa_t\Big(\frac{1}{J}\f_f\f\Big),J\Tt_f\right)_{L^2(\Omega_f)}\,dt + \int_0^T\left(\pa_t J\,\Tt_f,\frac{1}{J}\f_f\f\right)_{L^2(\Omega_f)}\,dt \\
        &= \int_0^T\left(\Tt_f,\f_f\right)_{L^2(\Omega_f)}\pa_t\f\,dt.
    \end{align*}
    This shows that $\Tt_f$ has the generalized time derivative $\pa_t\Tt_f=\frac{1}{J}\,\pa_t(J\Tt_f)-\frac{1}{J}\,\pa_t J\,\Tt_f$, which is in $L^2(0,T;H^1(\Omega_f)')$, see e.g.~\cite[Proposition 23.20]{zeidler_nonlinear_1990}.
    \qed
\end{proo}

Now, we can state the existence and uniqueness result for the advection-diffusion equation.

\begin{prop} \label{prop:theta_solution_estimates}
    There exists a constant $C>0$, only depending on the data from Assumption \ref{assu:data}, such that for all $\widetilde{\Tt}_s\in L^2(0,T;L^2(\Omega_s))$, the advection-diffusion equation \eqref{eq:transport_weak_ref} with coefficients given by \eqref{eq:ode_ref}, \eqref{eq:T_representation}, \eqref{eq:stokes_weak_ref}, \eqref{eq:R_representation}, \eqref{eq:S_representation}, and \eqref{eq:coefficients} admits a unique solution
    \begin{equation} \label{eq:theta_solution_spaces}
        \begin{alignedat}{2}
            &\Tt_f\in L^2(0,T;H^1(\Omega_f)),\qquad& 
            &\pa_t(J\Tt_f)\in L^2(0,T;H^1(\Omega_f)'),\\
            &\Tt_s\in L^2(0,T;H^1_{\Gamma_D}(\Omega_s)),&
            &\pa_t(J\Tt_s)\in L^2(0,T;H^1_{\Gamma_D}(\Omega_s)'),
        \end{alignedat}
    \end{equation}
    satisfying
    \begin{equation} \label{eq:theta_estimates}
        \N{\Tt_f}_{L^2(0,T;H^1(\Omega_f))} + \N{\pa_t\Tt_f}_{L^2(0,T;H^1(\Omega_f)')} + \N{\Tt_s}_{L^2(0,T;H^1_{\Gamma_D}(\Omega_s))} + \N{\pa_t\Tt_s}_{L^2(0,T;H^1_{\Gamma_D}(\Omega_s)')}\leq C.
    \end{equation}
\end{prop}

We prove Proposition \ref{prop:theta_solution_estimates} using the Galerkin method, see, for example, \cite[Chapter 23]{zeidler_nonlinear_1990}. Here, two specific difficulties arise from our problem: First, the domain is decomposed into $\Omega_f$ and $\Omega_s$, so we have the unknown tupel $(\Tt_f,\Tt_s)$ from the product space $H^1(\Omega_f)\times H^1_{\Gamma_D}(\Omega_s)$. Second, we need to treat the time and space dependent coefficients resulting from the nonlinear couplings in the system. In particular, the coefficients depend on the deformation gradient, and the time derivative in the Galerkin equations applies to the products $J\Tt_f^m$ and $J\Tt_s^m$, respectively, see also \cite{mabuza_modeling_2016}. Moreover, the convective term involves the solution of the Stokes system. \par
To shorten the notation, we will write the weak formulation \eqref{eq:transport_weak_ref} in terms of the time-dependent bilinear mapping
\begin{equation} \label{eq:addiff_a_def}
    \begin{aligned}
        &a:(H^1(\Omega_f)\times H^1_{\Gamma_D}(\Omega_s))\times(H^1(\Omega_f)\times H^1_{\Gamma_D}(\Omega_s))\times [0,T]\to\R, \\
        &a((\Tt_f,\Tt_s),(\f_f,\f_s),\cdot) := \int_{\Omega_f} (K^F_f\G\Tt_f \,-\,\Tt_f Aw)\cdot\G\f_f\,dx + \int_{\Omega_s} K^F_s\G\Tt_s\cdot\G\f_s\,dx \\
        &\hspace{4.5cm} +\int_{\Gamma_f}\Tt_f(w\cdot n)^+\f_f J\,d\sigma + \int_\Sigma \alpha(\Tt_f-\Tt_s)(\f_f-\f_s)J|F^{-T}n|\,d\sigma,
    \end{aligned}
\end{equation}
and the time-dependent linear mapping
\begin{equation} \label{eq:addiff_b_def}
    \begin{aligned}
        b:(H^1(\Omega_f)\times H^1_{\Gamma_D}(\Omega_s))\times [0,T]\to\R,~b((\f_f,\f_s),\cdot) := -\int_{\Gamma_f}f_{in} (w\cdot n)^-\f_f J\,d\sigma.
    \end{aligned}
\end{equation}
The weak formulation \eqref{eq:transport_weak_ref} now reads: Find $(\Tt_f,\Tt_s)$ in the spaces \eqref{eq:theta_solution_spaces} such that for all $(\f_f,\f_s)\in H^1(\Omega_f)\times H^1_{\Gamma_D}(\Omega_s)$ and a.e.~$t\in(0,T)$,
\begin{align*}
    \langle\pa_t(J(t)\Tt_f(t)),\f_f\rangle_{H^1(\Omega_f)',H^1(\Omega_f)} + \langle\pa_t(J(t)\Tt_s(t)),\f_s\rangle_{H^1_{\Gamma_D}(\Omega_s)',H^1_{\Gamma_D}(\Omega_s)}~+& \\
    +~a((\Tt_f(t),\Tt_s(t)),(\f_f,\f_s),t) &= b((\f_f,\f_s),t), \\
    (\Tt_f(0),\Tt_s(0))&=(\Tt_f^0,\Tt_s^0).
\end{align*}

\begin{rema}
    In this section, the unknown is the tuple $(\Tt_f,\Tt_s)$, so the (bi)linear mappings \eqref{eq:addiff_a_def} and \eqref{eq:addiff_b_def} are defined on the product space $H^1(\Omega_f)\times H^1_{\Gamma_D}(\Omega_s)$. 
    Accordingly, the test functions $(\f_f,\f_s)$ will be from the same product space. In particular, we may take test functions of the form $(\f_f,0)$ or $(0,\f_s)$ to show that a functional is in $L^2(0,T;H^1(\Omega_f)')$ or $L^2(0,T;H^1_{\Gamma_D}(\Omega_s)')$.
\end{rema}

\begin{proo}
    First, we prove the key properties of the (bi)linear mappings $a$ and $b$. Choose any $\widetilde{\Tt}_s\in L^2(0,T;L^2(\Omega_s))$ and let the coefficients in \eqref{eq:addiff_a_def} and \eqref{eq:addiff_b_def} be given by Propositions \ref{prop:S_solution_estimates} and \ref{prop:w_solution_estimates}. Then, for all $(\xi_f,\xi_s),(\f_f,\f_s)\in H^1(\Omega_f)\times H^1_{\Gamma_D}(\Omega_s)$, the function 
    \begin{equation*}
        a((\xi_f,\xi_s),(\f_f,\f_s),\cdot): [0,T]\to\R,~t\mapsto a((\xi_f,\xi_s),(\f_f,\f_s),t)
    \end{equation*}
    is continuous. Moreover, $a$ and $b$ are bounded. Indeed, for all $(\xi_f,\xi_s),(\f_f,\f_s)\in H^1(\Omega_f)\times H^1_{\Gamma_D}(\Omega_s)$ and all $t\in[0,T]$, we have
    \begin{align}
        &|a((\xi_f,\xi_s),(\f_f,\f_s),t)|\notag\\
        &\leq \N{K_f^{F(t)}}_{C(\ov{\Omega_f})}\N{\G\xi_f}_{L^2(\Omega_f)}\N{\G\f_f}_{L^2(\Omega_f)} + \N{K_s^{F(t)}}_{C(\ov{\Omega_s})}\N{\G\xi_s}_{L^2(\Omega_s)}\N{\G\f_s}_{L^2(\Omega_s)} \notag\\
        &\quad + \N{\xi_f}_{L^4(\Omega_f)}\N{A(t)}_{C(\ov{\Omega_f})}\N{w(t)}_{L^4(\Omega_f)}\N{\G\f_f}_{L^2(\Omega_f)} \notag\\
        &\quad + \N{\xi_f}_{L^4(\Gamma_f)}\N{w(t)}_{L^4(\Gamma_f)}\N{\f_f}_{L^2(\Gamma_f)}\N{J(t)}_{C(\Gamma_f)} \notag\\
        &\quad + \alpha\left(\N{\xi_f}_{L^2(\Sigma)} + \N{\xi_s}_{L^2(\Sigma)}\right) \left(\N{\f_f}_{L^2(\Sigma)} + \N{\f_s}_{L^2(\Sigma)}\right)\N{A(t)}_{C(\Sigma)} \notag\\
        &\leq \widetilde{C}\left(\N{K_f^{F(t)}}_{C(\ov{\Omega_f})} + \N{K_s^{F(t)}}_{C(\ov{\Omega_s})} + \N{A(t)}_{C(\ov{\Omega_f})}\N{w(t)}_{H^1(\Omega_f)} +\right. \notag\\
        &\hspace{1cm} + \left.\N{J(t)}_{C(\ov{\Omega_f})}\N{w(t)}_{H^1(\Omega_f)} + \N{A(t)}_{C(\Sigma)}\right) \notag \\
        &\quad \cdot \left(\N{\xi_f}_{H^1(\Omega_f)} + \N{\xi_s}_{H^1(\Omega_s)}\right) \left(\N{\f_f}_{H^1(\Omega_f)} + \N{\f_s}_{H^1(\Omega_s)}\right) \notag \\
        &\leq C\left(\N{\xi_f}_{H^1(\Omega_f)} + \N{\xi_s}_{H^1(\Omega_s)}\right) \left(\N{\f_f}_{H^1(\Omega_f)} + \N{\f_s}_{H^1(\Omega_s)}\right) \label{eq:addiff_a_continuity}
    \end{align}
    and, by the identification $L^2(0,T;H^1(\Omega_f)')=L^2(0,T;H^1(\Omega_f))'$,
    \begin{align}
        &\N{b((\cdot,0),\cdot)}_{L^2(0,T;H^1(\Omega_f)')} \notag\\
        &= \sup_{\f_f\in L^2(0,T;H^1(\Omega_f))} \N{\f_f}^{-1}_{L^2(0,T;H^1(\Omega_f))}\left|\int_0^T b((\f_f(t),0),t)\,dt\right| \notag\\
        &\leq \sup_{\f_f} \N{\f_f}^{-1}_{L^2(0,T;H^1(\Omega_f))}\int_0^T \N{f_{in}(t)}_{L^2(\Gamma_f)}\N{w(t)}_{L^4(\Gamma_f)}\N{\f_f(t)}_{L^4(\Gamma_f)}\N{J(t)}_{C(\Gamma_f)}\,dt \notag\\
        &\leq \widetilde{C}\sup_{\f_f} \N{\f_f}^{-1}_{L^2(0,T;H^1(\Omega_f))}\int_0^T \N{f_{in}(t)}_{L^2(\Gamma_f)}\N{w(t)}_{H^1(\Omega_f)}\N{\f_f(t)}_{H^1(\Omega_f)}\N{J(t)}_{C(\ov{\Omega_f})}\,dt \notag\\
        &\leq \widetilde{C}\N{f_{in}}_{L^2(0,T;L^2(\Gamma_f))}\N{w}_{C([0,T],H^1(\Omega_f))}\N{J}_{C([0,T]\times\ov{\Omega_f})} \notag\\
        &\leq C, \label{eq:addiff_b_continuity}
    \end{align}
    where in both estimates, $\widetilde{C}>0$ arises from the embeddings and $C>0$ depends also on the data from Assumption \ref{assu:data}. Finally, we derive the G\r{a}rding inequality
    \begin{equation} \label{eq:garding_inequality}
        \begin{aligned}
            &a((\xi_f,\xi_s),(\xi_f,\xi_s),\cdot) + \frac{1}{2}\int_{\Omega_f} \pa_t J|\xi_f|^2\,dx + \frac{1}{2}\int_{\Omega_s} \pa_t J|\xi_s|^2\,dx \\
            &\geq (c_f^K + c_s^K)\left(\N{\xi_f}_{H^1(\Omega_f)}^2 + \N{\xi_s}_{H^1(\Omega_s)}^2\right) - (c_f^K + c_s^K + C) \left(\N{\xi_f}_{L^2(\Omega_f)}^2 + \N{\xi_s}_{L^2(\Omega_s)}^2\right).
        \end{aligned}
    \end{equation}
    In $[0,T]$ and for all $(\xi_f,\xi_s)\in H^1(\Omega_f)\times H^1_{\Gamma_D}(\Omega_s)$, we have
    \begin{align*}
        &a((\xi_f,\xi_s),(\xi_f,\xi_s),\cdot) + \frac{1}{2}\int_{\Omega_f} \pa_t J|\xi_f|^2\,dx + \frac{1}{2}\int_{\Omega_s} \pa_t J|\xi_s|^2\,dx \\
        &\geq c_f^K\N{\G\xi_f}_{L^2(\Omega_f)}^2 + c_s^K\N{\G\xi_s}_{L^2(\Omega_s)}^2 + \frac{\alpha R_1}{R_0}\N{\xi_f-\xi_s}_{L^2(\Sigma)}^2 \\
        &\quad - \int_{\Omega_f} \xi_f Aw\cdot\G\xi_f\,dx + \int_{\Gamma_f} |\xi_f|^2(w\cdot n)^+ J\,d\sigma + \frac{1}{2}\int_{\Omega_f} \pa_t J|\xi_f|^2\,dx + \frac{1}{2}\int_{\Omega_s} \pa_t J|\xi_s|^2\,dx,
    \end{align*}
    where we used the constants from Proposition \ref{prop:S_solution_estimates}. For the last four terms, integration by parts yields:
    \begin{equation} \label{eq:convective_term_estimates}
        \begin{aligned}
            &\frac{1}{2}\int_{\Omega_f} Aw\cdot\G|\xi_f|^2\,dx - \int_{\Gamma_f} |\xi_f|^2(w\cdot n)^+ J\,d\sigma - \frac{1}{2}\int_{\Omega_f} \pa_t J|\xi_f|^2\,dx - \frac{1}{2}\int_{\Omega_s} \pa_t J|\xi_s|^2\,dx \\
            &= \frac{1}{2}\int_{\Gamma_f} |\xi_f|^2 Aw\cdot n\,d\sigma + \underbrace{\frac{1}{2}\int_{\Sigma} |\xi_f|^2 Aw\cdot n\,d\sigma}_{=0\text{ as }w=0\text{ on }\Sigma} - \frac{1}{2}\int_{\Omega_f} \underbrace{\Div(Aw)}_{=-\Div(Av_b),\text{ see \eqref{eq:stokes_weak_ref}.}}|\xi_f|^2\,dx \\
            &\quad - \int_{\Gamma_f} |\xi_f|^2(w\cdot n)^+ J\,d\sigma - \frac{1}{2}\int_{\Omega_f} \pa_t J|\xi_f|^2\,dx - \frac{1}{2}\int_{\Omega_s} \pa_t J|\xi_s|^2\,dx \\
            &= \frac{1}{2}\int_{\Gamma_f} |\xi_f|^2 (w\cdot \underbrace{F^{-T}n}_{=n\text{ on }\Gamma_f})J\,d\sigma + \frac{1}{2}\int_{\Omega_f} \underbrace{\Div(Av_b)}_{=\pa_t J}|\xi_f|^2\,dx \\
            &\quad - \int_{\Gamma_f} |\xi_f|^2(w\cdot n)^+ J\,d\sigma - \frac{1}{2}\int_{\Omega_f} \pa_t J|\xi_f|^2\,dx - \frac{1}{2}\int_{\Omega_s} \pa_t J|\xi_s|^2\,dx \\
            &= \underbrace{\frac{1}{2}\int_{\Gamma_f} |\xi_f|^2 (w\cdot n)^- J\,d\sigma}_{\leq 0} \underbrace{- \frac{1}{2}\int_{\Gamma_f} |\xi_f|^2(w\cdot n)^+ J\,d\sigma}_{\leq 0} \\
            &\quad + \frac{1}{2}\int_{\Omega_f} \pa_t J|\xi_f|^2\,dx - \frac{1}{2}\int_{\Omega_f} \pa_t J|\xi_f|^2\,dx - \frac{1}{2}\int_{\Omega_s} \pa_t J|\xi_s|^2\,dx \\
            &\leq -\frac{1}{2}\int_{\Omega_s} \pa_t J|\xi_s|^2\,dx \\
            &\leq \N{\pa_t J}_{C(\ov{\Omega_s})} \N{\xi_s}_{L^2(\Omega_s)}^2 \\
            &\leq C \N{\xi_s}_{L^2(\Omega_s)}^2,
        \end{aligned}
    \end{equation}
    with some uniform constant $C>0$ from Proposition \ref{prop:S_solution_estimates}. Note that $\Div(Av_b)=\pa_t J$ follows from Euler's expansion formula, see e.g.~\cite[Satz 5.2]{eck_mathematische_2017}, and the Piola identity \eqref{eq:piola_identity}. Moreover, $F^{-T}n=n$ on $\Gamma_f$ is due to \eqref{eq:F_representation}. In summary, we obtain \eqref{eq:garding_inequality}.\par
    Now, let $(v_k)_k$ and $(w_k)_k$ be two Schauder bases for $H^1(\Omega_f)$ and $H^1_{\Gamma_D}(\Omega_s)$, respectively, and let $(\Tt_f^{0,m})_m\subset L^2(\Omega_f)$ and $(\Tt_s^{0,m})_m\subset L^2(\Omega_s)$ be such that
    \begin{alignat*}{3}
        &\Tt_f^{0,m}\to\Tt_f^0\text{ in }L^2(\Omega_f), & &\Tt_s^{0,m}\to\Tt_s^0\text{ in }L^2(\Omega_s) & &\text{as }m\to\infty, \\
        &\Tt_f^{0,m}\in\mathrm{span}\{v_1,\dots,v_m\}, & \qquad &\Tt_s^{0,m}\in\mathrm{span}\{w_1,\dots,w_m\} & \qquad &\text{for all }m\in\NN.
    \end{alignat*}
    As $(v_1,\dots,v_m)$ and $(w_1,\dots,w_m)$ are linearly independent (also in $L^2(\Omega_f)$ and $L^2(\Omega_s)$) for all $m\in\NN$, there exist unique $\alpha_i^m,\beta_i^m\in\R$ (for $m\in\NN$ and $1\leq i\leq m$) such that
    \begin{equation*}
        \Tt_f^{0,m} = \sum_{i=1}^m \alpha_i^m v_i \qquad\text{and}\qquad \Tt_s^{0,m} = \sum_{i=1}^m \beta_i^m w_i.
    \end{equation*}
    For fixed $m\in\NN$, we define the Galerkin equations by
    \begin{equation} \label{eq:galerkin_equations}
        \begin{aligned}
            \frac{d}{dt}\left(\sum_{i=1}^m a_i^m(t) \int_{\Omega_f} J(t)v_iv_k\,dx\right) + \sum_{i=1}^m a((a_i^m(t)v_i,b_i^m(t)w_i),(v_k,0),t) &= b((v_k,0),t), \\
            \frac{d}{dt}\left(\sum_{i=1}^m b_i^m(t) \int_{\Omega_s} J(t)w_iw_k\,dx\right) + \sum_{i=1}^m a((a_i^m(t)v_i,b_i^m(t)w_i),(0,w_k),t) &= b((0,w_k),t), \\
            a_k^m(0) &= \alpha_k^m, \\
            b_k^m(0) &= \beta_k^m,
        \end{aligned}
    \end{equation}
    for all $t\in[0,T]$ and all $k=1,\dots,m$. Note that \eqref{eq:galerkin_equations} can be formally derived from \eqref{eq:transport_weak_ref} by replacing $\Tt_f$ by $\Tt_f^m(t):=\sum_{i=1}^m a_i^m(t)v_i$, $\Tt_s$ by $\Tt_s^m(t):=\sum_{i=1}^m b_i^m(t)w_i$, $(\Tt_f^0,\Tt_s^0)$ by $(\Tt_f^{0,m},\Tt_s^{0,m})$, and $(\f_f,\f_s)$ by $(v_k,0)$ or $(0,w_k)$ for $k=1,\dots,m$. Using the biliearity of $a$, the system \eqref{eq:galerkin_equations} can be written as
    \begin{equation} \label{eq:galerkin_matrix_form}
        \begin{aligned}
            \frac{d}{dt}(L^m(t)c^m(t)) + M^m(t)c^m(t) &= N^m(t), \\
            c^m(0) &= \gamma^m,
        \end{aligned}
    \end{equation}
    for all $t\in[0,T]$, where:
    \begin{align*}
        L^m(t) &= \begin{pmatrix}
            L_v^m(t) & 0 \\
            0 & L_w^m(t)
        \end{pmatrix}
        = \begin{pmatrix}
            \left(\int_{\Omega_f}J(t)v_i v_k\,dx\right)_{k,i=1,\dots,m} & 0 \\
            0 & \left(\int_{\Omega_s}J(t)w_i w_k\,dx\right)_{k,i=1,\dots,m}
        \end{pmatrix}, \\
        M^m(t) &= \begin{pmatrix}
            (a((v_i,0),(v_k,0),t))_{k,i=1,\dots,m} & (a((0,w_i),(v_k,0),t))_{k,i=1,\dots,m} \\
            (a((v_i,0),(0,w_k),t))_{k,i=1,\dots,m} & (a((0,w_i),(0,w_k),t))_{k,i=1,\dots,m}
        \end{pmatrix}, \\
        N^m(t) &= \begin{pmatrix}
            (b((v_k,0),t))_{k=1,\dots,m} \\
            (b((0,w_k),t))_{k=1,\dots,m}
        \end{pmatrix}, \\
        c^m(t) &= \begin{pmatrix}
            (a_i^m(t))_{i=1\dots,m} \\
            (b_i^m(t))_{i=1\dots,m}
        \end{pmatrix}, \\
        \gamma^m(t) &= \begin{pmatrix}
            (\alpha_i^m(t))_{i=1\dots,m} \\
            (\beta_i^m(t))_{i=1\dots,m}
        \end{pmatrix}.
    \end{align*}
    Clearly, $L^m(t)$ is symmetric. To prove that it is positive definite for all $t\in[0,T]$, take any $\xi\in\R^m\setminus\{0\}$ and compute:
    \begin{equation*}
        \xi^T L_v^m(t)\xi = \sum_{i,k=1}^m\int_{\Omega_f} J(t)v_i v_k \xi_i \xi_k\,dx = \int_{\Omega_f} J(t)\left(\sum_{i=1}^m \xi_i v_i\right)^2\,dx \geq \eta\N{\sum_{i=1}^m \xi_i v_i}_{L^2(\Omega_f)}^2 \geq 0
    \end{equation*}
    If we had $\xi^T L_v^m(t)\xi = 0$, then $\sum_{i=1}^m \xi_i v_i = 0$ a.e.~in $\Omega_f$. But this would imply $\xi = 0$ as $(v_i)_i$ is a basis for $H^1(\Omega_f)$. Consequently, we have $\xi^T L_v^m(t)\xi > 0$, and in the same way, we get $\xi^T L_w^m(t)\xi > 0$, i.e.~$L^m(t)$ is positive definite for all $t\in[0,T]$. Thus, we can substitute $c^m(t) = (L^m)^{-1}(t)d^m(t)$ in \eqref{eq:galerkin_matrix_form} and solve the system
    \begin{equation} \label{eq:galerkin_substituted_form}
        \begin{aligned}
            \frac{d}{dt}d^m(t) + M^m(t)(L^m)^{-1}(t)d^m(t) &= N^m(t), \\
            d^m(t) &= L^m(0)\gamma^m,
        \end{aligned}
    \end{equation}
    for all $t\in[0,T]$. Due to Propositions \ref{prop:S_solution_estimates} and \ref{prop:w_solution_estimates}, $L^m$ and $M^m$ are continuous in time in $[0,T]$, so also $(L^m)^{-1}$ is continuous in time in $[0,T]$. Moreover, Assumption \ref{assu:data} (A5) guarantees that $N^m\in L^2(0,T)$. As the first line of \eqref{eq:galerkin_substituted_form} is a linear system, it satisfies the Carathéodory conditions and thus, by Carathéodory's existence theorem, the initial value problem \eqref{eq:galerkin_substituted_form} admits a unique absolutely continuous solution, see, e.g., \cite[Section I.5]{hale_ordinary_1969}. Moreover, $L^m$ and its inverse are absolutely continuous (even continuously differentiable) in $[0,T]$, so we obtain an absolutely continuous solution $c^m(t) = (L^m)^{-1}(t)d^m(t)$ of \eqref{eq:galerkin_matrix_form}. To derive uniform estimates for the $m$-th Galerkin solution $(\Tt_f^m,\Tt_s^m)$, multiply the Galerkin equations \eqref{eq:galerkin_equations} by $a_k^m(t)$ resp. $b_k^m(t)$, sum both equations and sum over $k=1,\dots,m$ to obtain:
    \begin{equation} \label{eq:addiff_tested_with_galerkin_solution}
        \begin{aligned}
            \int_{\Omega_f} \pa_t(J(t)\Tt_f^m(t))\,\Tt_f^m(t)\,dx + \int_{\Omega_s} \pa_t(J(t)\Tt_s^m(t))\,\Tt_s^m(t)\,dx~+&\\
            +~a((\Tt_f^m(t),\Tt_s^m(t)),(\Tt_f^m(t),\Tt_s^m(t)),t) &= b((\Tt_f^m(t),\Tt_s^m(t)),t).
        \end{aligned}
    \end{equation}
    To conclude the uniform estimates, use the product rule for the time derivative and apply the G\r{a}rding inequality \eqref{eq:garding_inequality}, the estimate \eqref{eq:addiff_b_continuity} for $b$, and the Gronwall inequality. \par
    For the rest of the proof, we can follow the standard procedure, see, for example, \cite[Section 23.9]{zeidler_nonlinear_1990}. After passing to the limit $m\to\infty$, we obtain limit functions $\Tt_f\in L^2(0,T;H^1(\Omega_f))$ and $\Tt_s\in L^2(0,T;H^1_{\Gamma_D}(\Omega_s))$ for which the products $J\Tt_f$ and $J\Tt_s$ have generalized time derivatives in $L^2(0,T;H^1(\Omega_f)')$ and $L^2(0,T;H^1_{\Gamma_D}(\Omega_s)')$. Thus $(\Tt_f,\Tt_s)$ is a solution of the advection-diffusion equation \eqref{eq:transport_weak_ref}. The uniform estimate \eqref{eq:theta_estimates} for $\pa_t\Tt_f$ and $\pa_t\Tt_s$ relies on the product rule \eqref{eq:product_rule_J} and the estimates \eqref{eq:addiff_a_continuity} and \eqref{eq:addiff_b_continuity}. Moreover, uniqueness follows from the product rule \eqref{eq:product_rule_sqrt_J}, the G\r{a}rding inequality \eqref{eq:garding_inequality}, and the Gronwall inequality.
    \qed
\end{proo}

\subsection{Existence of solutions for the fully coupled system} \label{subsec:fixed_point_method}

Using Proposition \ref{prop:theta_solution_estimates}, we now define the solution operator $\mathcal{F}$ of the linearized system by \eqref{eq:solution_operator}. In order to apply Schaefer's fixed-point theorem, we need the following lemma.

\begin{lemm} \label{lemm:F_continuity}
    The mapping $\mathcal{F}$ is compact and continuous.
\end{lemm}

\begin{proo}
    To show compactness, let $(\widetilde{\Tt}_s^k)_k$ be a bounded sequence in $L^2(0,T;L^2(\Omega_s))$. By Proposition \ref{prop:theta_solution_estimates}, there exists a constant $C>0$ such that
    \begin{equation*}
        \N{\mathcal{F}(\widetilde{\Tt}_s^k)}_{L^2(0,T;H^1_{\Gamma_D}(\Omega_s))} + \N{\pa_t\mathcal{F}(\widetilde{\Tt}_s^k)}_{L^2(0,T;H^1_{\Gamma_D}(\Omega_s)')}\leq C
    \end{equation*}
    for all $k\in\NN$. By the Aubin--Lions lemma, $(\mathcal{F}(\widetilde{\Tt}_s^k))_k$ is relatively compact in $L^2(0,T;L^2(\Omega_s))$, see, e.g., \cite[Corollary 4]{simon_compact_1986}, so $\mathcal{F}$ is compact. \par
    The continuity requires some more work. Let $(\widetilde{\Tt}_s^k)_k\subset L^2(0,T;L^2(\Omega_s))$ converge to $\widetilde{\Tt}_s$ in $L^2(0,T;L^2(\Omega_s))$. We denote by
    \begin{equation*}
        (c^k,w^k,q^k,\Tt_f^k,\Tt_s^k)\qquad\text{and}\qquad(c,w,q,\Tt_f,\Tt_s)
    \end{equation*}
    the corresponding unique solutions given by Propositions \ref{prop:c_solution_estimates}, \ref{prop:w_solution_estimates}, and \ref{prop:theta_solution_estimates}. Moreover, we write
    \begin{equation*}
        (R^k,S_k,v_b^k,F_k,J_k,A_k)\qquad\text{and}\qquad(R,S,v_b,F,J,A)
    \end{equation*}
    for the corresponding radius, deformation and coefficients given by \eqref{eq:R_representation}, \eqref{eq:S_representation}, and \eqref{eq:coefficients}. To prove the continuity of $\mathcal{F}$, we need to show $\Tt_s^k\to\Tt_s$ in $L^2(0,T;L^2(\Omega_s))$. \par
    Using Young's convolution inequality and the representation of $\mathcal{T}$ from Lemma \ref{lemm:T_properties}, we obtain:
    \begin{align*}
        \N{\mathcal{T}(\widetilde{\Tt}_s^k)-\mathcal{T}(\widetilde{\Tt}_s)}_{C([0,T])} &= \N{K_\gamma\star\left(\mathcal{T}_1(\widetilde{\Tt}_s^k)-\mathcal{T}_1(\widetilde{\Tt}_s)\right)}_{C([0,T])} \\
        &\leq C \N{\mathcal{T}_1(\widetilde{\Tt}_s^k)-\mathcal{T}_1(\widetilde{\Tt}_s)}_{L^1(-\gamma,T)} = C \N{\mathcal{T}_1(\widetilde{\Tt}_s^k)-\mathcal{T}_1(\widetilde{\Tt}_s)}_{L^1(0,T)} \\
        &\leq \widetilde{C} \N{\widetilde{\Tt}_s^k-\widetilde{\Tt}_s}_{L^2(0,T;L^2(\Omega_s))} \to 0\qquad\text{as }k\to\infty.
    \end{align*}
    From Assumption \ref{assu:data} (A2), it follows that $G$, $\pa_{x_1}G$, and $\pa_{x_1}^2 G$ are Lipschitz continuous. Using the representation \eqref{eq:c_solution} from Proposition \ref{prop:c_solution_estimates}, we obtain:
    \begin{align*}
        &\sup_{(t,x_1)\in[0,T]\times[0,L]} |c^k(t,x_1) - c(t,x_1)| \\
        &=\sup_{(t,x_1)\in[0,T]\times[0,L]}\EN{e^{-kt}\int_0^t e^{ks}\big(G(x_1,\mathcal{T}(\widetilde{\Tt}_s^k)(s)) - G(x_1,\mathcal{T}(\widetilde{\Tt}_s)(s))\big)\,ds} \\
        &\leq C\sup_{(t,x_1)\in[0,T]\times[0,L]}\EN{G(x_1,\mathcal{T}(\widetilde{\Tt}_s^k)(t)) - G(x_1,\mathcal{T}(\widetilde{\Tt}_s)(t))} \\
        &\leq\widetilde{C}\sup_{(t,x_1)\in[0,T]\times[0,L]}\EN{\mathcal{T}(\widetilde{\Tt}_s^k)(t) - \mathcal{T}(\widetilde{\Tt}_s)(t))} \to 0\qquad\text{as }k\to\infty,
    \end{align*}
    and with \eqref{eq:c_time_derivative}:
    \begin{align*}
        &\sup_{(t,x_1)\in[0,T]\times[0,L]} |\pa_t c^k(t,x_1) - \pa_t c(t,x_1)| \\
        &=\sup_{(t,x_1)\in[0,T]\times[0,L]}\Big|-k\,e^{-kt}\int_0^t e^{ks}\big(G(x_1,\mathcal{T}(\widetilde{\Tt}_s^k)(s)) - G(x_1,\mathcal{T}(\widetilde{\Tt}_s)(s))\big)\,ds \\
        &\hspace{3.2cm}+ G(x_1,\mathcal{T}(\widetilde{\Tt}_s^k)(s)) - G(x_1,\mathcal{T}(\widetilde{\Tt}_s)(s))\Big| \\
        &\leq (C+1)\sup_{(t,x_1)\in[0,T]\times[0,L]}\EN{G(x_1,\mathcal{T}(\widetilde{\Tt}_s^k)(t)) - G(x_1,\mathcal{T}(\widetilde{\Tt}_s)(t))} \\
        &\leq\widetilde{C}\sup_{(t,x_1)\in[0,T]\times[0,L]}\EN{\mathcal{T}(\widetilde{\Tt}_s^k)(t) - \mathcal{T}(\widetilde{\Tt}_s)(t))} \to 0\qquad\text{as }k\to\infty.
    \end{align*}
    Replacing $G$ by $\pa_{x_1}G$ or $\pa_{x_1}^2 G$, it follows
    \begin{align*}
        &\sup_{(t,x_1)\in[0,T]\times[0,L]} |\pa_{x_1} c^k(t,x_1) - \pa_{x_1} c(t,x_1)|\to 0\qquad\text{as }k\to\infty, \\
        &\sup_{(t,x_1)\in[0,T]\times[0,L]} |\pa_{x_1}^2 c^k(t,x_1) - \pa_{x_1}^2 c(t,x_1)|\to 0\qquad\text{as }k\to\infty, \\
        &\sup_{(t,x_1)\in[0,T]\times[0,L]} |\pa_t\pa_{x_1} c^k(t,x_1) - \pa_t\pa_{x_1} c(t,x_1)|\to 0\qquad\text{as }k\to\infty,
    \end{align*}
    so altogether, $c^k\to c$ in $C_1^2([0,T]\times[0,L])$ as $k\to\infty$. \par
    By $R^k=H\circ c^k$ and $R=H\circ c$, see \eqref{eq:R_representation}, the convergence of $(c^k)_k$ implies $R^k\to R$ in $C_1^2([0,T]\times[0,L])$ as $k\to\infty$ due to the boundedness and Lipschitz continuity of $H$, $\pa_y H$, and $\pa_y^2 H$, see Assumption \ref{assu:data} (A3). \par
    In the same way, we get $S^k\to S$ in $C_1^2([0,T]\times\ov{\Omega})$ as $k\to\infty$ due to the boundedness and Lipschitz continuity of $\pa_R^{\beta_1}\pa_r^{\beta_2}\rho$ in $[R_1,R_2]\times[0,\frac{1}{2}]$ for all $\beta\in\NN_0^2$ with $0\leq |\beta|\leq 2$, see Remark \ref{rema:S_bijectivity} (2). Note that it suffices to prove the uniform convergence of $(S^k)_k$ and its derivatives in $[0,T]\times\ov{Z}$ instead of $[0,T]\times\ov{\Omega}$, see \eqref{eq:Z_definition} and Remark \ref{rema:S_bijectivity} (3). Moreover, this implies $v_b^k\to v_b$, $F_k\to F$, $J_k\to J$, and $A_k\to A$ in $C([0,T]\times\ov{\Omega})$ as $k\to\infty$. \par
    In order to prove the convergence of $(w^k,q^k)_k$, let
    \begin{equation*}
        (a^k,b^k,f^k,g^k)\qquad\text{and}\qquad(a,b,f,g)
    \end{equation*}
    be the (bi)linear forms \eqref{eq:stokes_bilinear_forms} corresponding to $\widetilde{\Tt}_s^k$ and $\widetilde{\Tt}_s$. By Proposition \ref{prop:w_solution_estimates}, there exist constants $\alpha_0,\beta_0,C>0$ such that \eqref{eq:abfg_solution} holds for $(a,b,f,g)$ and $(a^k,b^k,f^k,g^k)$ for all $k\in\NN$. The following estimates resemble those in the proof of Proposition \ref{prop:w_solution_estimates} where we have proven that $a$, $b$, $f$, and $g$ are continuous in time, so some steps will be omitted here. For all $t\in[0,T]$ and all $u,v\in \V$, we have:
    \begin{align*}
        |a^k(t)(u,v)-a(t)(u,v)| =\dots&\leq \mu\widetilde{C}\N{u}_\V\N{v}_\V\left(\N{F_k^{-T}(t)-F^{-T}(t)}_{C(\ov{\Omega_f})}\N{A_k^T(t)}_{C(\ov{\Omega_f})}\right.\\
        &\hspace{3.3cm} + \left.\N{F^{-T}(t)}_{C(\ov{\Omega_f})}\N{A_k^T(t)-A^T(t)}_{C(\ov{\Omega_f})}\right),
    \end{align*}
    and due to $S^k\to S$ in $C_1^2([0,T]\times\ov{\Omega})$ as $k\to\infty$, we conclude
    \begin{align*}
        \sup_{t\in[0,T]}\N{a^k(t)-a(t)}_\Bil &\leq \mu\widetilde{C}\left(\N{F_k^{-T}-F^{-T}}_{C([0,T]\times\ov{\Omega_f})}\N{A_k^T}_{C([0,T]\times\ov{\Omega_f})}\right.\\
        &\hspace{1.3cm} + \left.\N{F^{-T}}_{C([0,T]\times\ov{\Omega_f})}\N{A_k^T-A^T}_{C([0,T]\times\ov{\Omega_f})}\right) \\
        &\to 0\qquad\text{as }k\to\infty.
    \end{align*}
    Similarly, for all $t\in[0,T]$, all $v\in \V$ and all $p\in \Q$,
    \begin{equation*}
        |-b^k(t)(v,p)+b(t)(v,p)| =\dots\leq \widetilde{C}\N{v}_\V\N{p}_\Q\N{A_k(t)-A(t)}_{C(\ov{\Omega_f})},
    \end{equation*}
    and thus
    \begin{equation*}
        \sup_{t\in[0,T]}\N{b^k(t)-b(t)}_\Bil\leq\widetilde{C}\N{A_k-A}_{C([0,T]\times\ov{\Omega_f})}\to 0\qquad\text{as }k\to\infty.
    \end{equation*}
    Using Poincaré's inequality, for all $t\in[0,T]$ and all $\phi\in \V$, we get:
    \begin{align*}
        &|-\DP{f^k(t),\phi}{\V} + \DP{f(t),\phi}{\V}|\\
        &=\Big|\Of \left(A_k^T(t)\G f_b(t) - A^T(t)\G f_b(t)\right)\cdot\phi\dx \\
        &\qquad+~a^k(t)(v_b^k(t),\phi) - a(t)(v_b^k(t),\phi) + a(t)(v_b^k(t)-v_b(t),\phi)\Big|\\
        &\leq \widetilde{C}\N{\phi}_\V\left(\N{A_k^T(t)-A^T(t)}_{C(\ov{\Omega_f})}\N{f_b(t)}_{H^1(\Omega_f)}\right.\\
        &\hspace{2.1cm} + \left.\mu\N{a^k(t)-a(t)}_\Bil\N{\G v_b^k(t)}_{C(\ov{\Omega_f})} +~\mu\N{a(t)}_\Bil\N{\G v_b^k(t) - \G v_b(t)}_{C(\ov{\Omega_f})}\right).
    \end{align*}
    and thus
    \begin{align*}
        \sup_{t\in[0,T]}\N{f^k(t)-f(t)}_{\V'}
        &\leq \widetilde{C}\Big(\N{A_k^T-A^T}_{C([0,T]\times\ov{\Omega_f})}\N{f_b}_{C([0,T],H^1(\Omega_f))}\\
        &\qquad + \mu\Big(\sup_{t\in[0,T]}\N{a^k(t)-a(t)}_\Bil\Big)\N{\G v_b^k}_{C([0,T]\times\ov{\Omega_f})} \\
        &\qquad + \mu\Big(\sup_{t\in[0,T]}\N{a(t)}_\Bil\Big)\N{\G v_b^k - \G v_b}_{C([0,T]\times\ov{\Omega_f})}\Big)\\
        &\to 0\qquad\text{as }k\to\infty.
    \end{align*}
    Furthermore, for all $\psi\in \Q$, we can estimate
    \begin{align*}
        |\DP{g^k(t),\psi}{\Q} - \DP{g(t),\psi}{\Q}|&=\Big|\Of \psi\,(A_k(t):\G v_b^k(t) - A(t):\G v_b(t))\dx\Big| \\
        &\leq\widetilde{C}\N{\psi}_\Q\N{A_k(t):\G v_b^k(t) - A(t):\G v_b(t)}_{C(\ov{\Omega_f})},
    \end{align*}
    so we end up with
    \begin{equation*}
        \N{g^k(t)-g(t)}_{\Q'}\leq\widetilde{C}\N{A_k:\G v_b^k - A:\G v_b}_{C([0,T]\times\ov{\Omega_f})} \to 0\qquad\text{as }k\to\infty.
    \end{equation*}
    Now, using \eqref{eq:abfg_solution} and the general Lipschitz estimate \cite[Lemma B.4]{gahn_rigorous_2025}, we can find a constant $L(\alpha_0,\beta_0,C)>0$, only depending on $\alpha_0$, $\beta_0$, and $C$, such that:
    \begin{align*}
        &\N{w^k-w}_{C([0,T],\V)} + \N{q^k-q}_{C([0,T],\Q)} \\
        &\leq 2\sup_{t\in[0,T]} \Big(\N{w^k(t)-w(t)}_\V + \N{q^k(t)-q(t)}_\Q\Big) \\
        &\leq 2\,L(\alpha_0,\beta_0,C)\sup_{t\in[0,T]} \Big(\N{a^k(t)-a(t)}_\Bil + \N{b^k(t)-b(t)}_\Bil \\
        &\hspace{4.2cm}+ \N{f^k(t)-f(t)}_{\V'} + \N{g^k(t)-g(t)}_{\Q'}\Big) \\
        &\to 0\qquad\text{as }k\to\infty.
    \end{align*} \par
    As a last step, we prove the convergence of $(\Tt_f^k,\Tt_s^k)_k$. In the following, we denote by
    \begin{equation*}
        a^k,b^k\qquad\text{and}\qquad a,b
    \end{equation*}
    the (bi)linear forms \eqref{eq:addiff_a_def}, \eqref{eq:addiff_b_def} corresponding to $\widetilde{\Tt}_s^k$ and $\widetilde{\Tt}_s$. For all $(\f_f,\f_s)\in H^1(\Omega_f)\times H^1_{\Gamma_D}(\Omega_s)$, all $\f\in C_c^\infty(0,T)$, and all $\zeta_f\in L^2(0,T;H^1(\Omega_f))$, the same embedding inequalities as in \eqref{eq:addiff_a_continuity} lead to the estimate
    \begin{align}
        &\int_0^T a^k((\zeta_f,0),(\f_f,\f_s),t)\f - a((\zeta_f,0),(\f_f,\f_s),t)\f\,dt \notag\\
        &= \int_0^T\int_{\Omega_f} ((K^{F_k}_f-K^{F}_f)\G\zeta_f \,-\,\zeta_f (A_k w^k-Aw))\cdot\G\f_f\,dx \notag\\
        &\quad +\int_{\Gamma_f}\zeta_f(J_k(w^k\cdot n)^+ -J(w\cdot n)^+)\f_f\,d\sigma + \int_\Sigma \alpha\,\zeta_f(\f_f-\f_s)(J_k|F_k^{-T}n|-J|F^{-T}n|)\,d\sigma\,\f\,dt \notag\\
        &\leq \widetilde{C}\left(\N{K^{F_k}_f-K^{F}_f}_{C([0,T]\times\ov{\Omega_f})} + \N{A_k}_{C([0,T]\times\ov{\Omega_f})}\N{w^k-w}_{C([0,T],H^1(\Omega_f))} +\right. \notag\\
        &\hspace{1cm} + \left.\N{A_k-A}_{C([0,T]\times\ov{\Omega_f})}\N{w}_{C([0,T],H^1(\Omega_f))} + \N{A_k-A}_{C([0,T]\times\Sigma)}\right) \label{eq:addiff_a_k_convergence}\\
        &\quad \cdot \N{\zeta_f}_{L^2(0,T;H^1(\Omega_f))}\left(\N{\f_f}_{H^1(\Omega_f)} + \N{\f_s}_{H^1(\Omega_s)}\right)\N{\f}_{C([0,T])}.\notag
    \end{align}
    Note that we have applied the reverse triangle inequality, the identity $F_k^{-T}n=F^{-T}n=n$ on $\Gamma_f$, see \eqref{eq:F_representation}, and $J_k,J\geq \eta>0$, see \eqref{eq:S_estimate_J}, in the following way:
    \begin{align*}
        &|J_k|F_k^{-T}n|-J|F^{-T}n||\leq |J_kF_k^{-T}-JF^{-T}|\underbrace{|n|}_{=1}\leq |A_k^T-A^T|\qquad\text{in }[0,T]\times\Sigma, \\
        &J_k(w^k\cdot n)^+ -J(w\cdot n)^+ = J_k\frac{|w^k\cdot n|+w^k\cdot n}{2} - J\frac{|w\cdot n|+w\cdot n}{2} \\
        &\leq |J_k w^k\cdot (F_k^{-T}n) - Jw\cdot (F^{-T}n)| \leq |A_kw^k - Aw|\qquad\text{in }[0,T]\times\Gamma_f.
    \end{align*}
    Consequently, we obtain
    \begin{equation*}
        \textstyle \N{\int_0^T a^k((\cdot,0),(\f_f,\f_s),t)\f - a((\cdot,0),(\f_f,\f_s),t)\f\,dt}_{L^2(0,T;H^1(\Omega_f)')}\to 0\qquad\text{as }k\to\infty,
    \end{equation*}
    and in a similar way,
    \begin{equation*}
        \textstyle \N{\int_0^T a^k((0,\cdot),(\f_f,\f_s),t)\f - a((0,\cdot),(\f_f,\f_s),t)\f\,dt}_{L^2(0,T;H_{\Gamma_D}^1(\Omega_s)')}\to 0\qquad\text{as }k\to\infty.
    \end{equation*}
    Moreover, the same trace estimates as in \eqref{eq:addiff_b_continuity} yield
    \begin{align}
        &\int_0^T b^k((\f_f,\f_s),t)\f - b((\f_f,\f_s),t)\f\,dt \notag\\
        &= -\int_0^T\int_{\Gamma_f}f_{in} (J_k(w^k\cdot n)^- - J(w\cdot n)^-)\f_f\,d\sigma\,\f\,dt \notag\\
        &\leq \widetilde{C}\N{f_{in}}_{L^2(0,T;L^2(\Gamma_f))}\left(\N{J_k}_{C([0,T]\times\ov{\Omega_f})}\N{w^k-w}_{C([0,T],H^1(\Omega_f))}\right. \label{eq:addiff_b_k_convergence}\\
        &\hspace{3.2cm} + \left.\N{J_k-J}_{C([0,T]\times\ov{\Omega_f})}\N{w}_{C([0,T],H^1(\Omega_f))}\right)\N{\f_f}_{H^1(\Omega_f)}\N{\f}_{C([0,T])} \notag\\
        &\to 0\qquad\text{as }k\to\infty.\notag
    \end{align}
    By Proposition \ref{prop:theta_solution_estimates}, there exists a constant $C>0$ such that \eqref{eq:theta_estimates} holds for $(\Tt_f^k,\Tt_s^k)$ for all $k\in\NN$. Thus, we can extract twice a subsequence, still denoted by $(\Tt_f^k,\Tt_s^k)_k$, such that:
    \begin{equation} \label{eq:theta_convergence}
        \begin{aligned}
            &\Tt_f^k\rightharpoonup\ov{\Tt}_f\qquad\text{in }L^2(0,T;H^1(\Omega_f)), \\
            &\Tt_s^k\rightharpoonup\ov{\Tt}_s\qquad\text{in }L^2(0,T;H_{\Gamma_D}^1(\Omega_s)), \\
            &\Tt_f^k\to\ov{\Tt}_f\qquad\text{in }L^2(0,T;L^2(\Omega_f)), \\
            &\Tt_s^k\to\ov{\Tt}_s\qquad\text{in }L^2(0,T;L^2(\Omega_s)),
        \end{aligned}
    \end{equation}
    as $k\to\infty$, where the strong convergence follows from Aubin-Lions lemma, see e.g.~\cite[Corollary 4]{simon_compact_1986}. As $(\Tt_f^k,\Tt_s^k)$ is the weak solution corresponding to $\widetilde{\Tt}_s^k$, we know that for all $(\f_f,\f_s)\in H^1(\Omega_f)\times H^1_{\Gamma_D}(\Omega_s)$ and all $\f\in C_c^\infty(0,T)$:
    \begin{align*}
        -\int_0^T\int_{\Omega_f} J_k\Tt_f^k\f_f\,\pa_t\f\,dx\,dt - \int_0^T\int_{\Omega_s} J_k\Tt_s^k\f_s\,\pa_t\f\,dx\,dt~+& \\
        + \int_0^T a^k((\Tt_f^k,0),(\f_f,\f_s),t)\f\,dt + \int_0^T a^k((0,\Tt_s^k),(\f_f,\f_s),t)\f\,dt &=\int_0^T b^k((\f_f,\f_s),t)\f\,dt.
    \end{align*}
    Remark that this weak formulation is equivalent to \eqref{eq:transport_weak_ref} by the definition of the distributional time derivative, see e.g.~\cite[Proposition 23.20]{zeidler_nonlinear_1990}. Now, using the convergence results from above, we can pass to the limit $k\to\infty$ in all terms and end up with
    \begin{align*}
        -\int_0^T\int_{\Omega_f} J\Tt_f\f_f\,\pa_t\f\,dx\,dt - \int_0^T\int_{\Omega_s} J\Tt_s\f_s\,\pa_t\f\,dx\,dt~+& \\
        + \int_0^T a((\ov{\Tt}_f,0),(\f_f,\f_s),t)\f\,dt + \int_0^T a((0,\ov{\Tt}_s),(\f_f,\f_s),t)\f\,dt &=\int_0^T b((\f_f,\f_s),t)\f\,dt,
    \end{align*}
    so $(\ov{\Tt}_f,\ov{\Tt}_s)$ is a weak solution corresponding to $\widetilde{\Tt}_s$. By uniqueness, it follows $(\ov{\Tt}_f,\ov{\Tt}_s) = (\Tt_f,\Tt_s)$, and thus, by \eqref{eq:theta_convergence}, $\Tt_s^k\to\Tt_s$ in $L^2(0,T;L^2(\Omega_s))$ as $k\to\infty$.
    \qed
\end{proo}

\begin{proo}[Existence for the fully coupled system \eqref{eq:solution_spaces}--\eqref{eq:coefficients}]
    First, let us consider the set
    \begin{equation*}
        \mathcal{E}:=\{\widetilde{\Tt}_s\in L^2(0,T;L^2(\Omega_s)): \widetilde{\Tt}_s = \lambda\mathcal{F}(\widetilde{\Tt}_s)\text{ for some }0\leq\lambda\leq 1\}.
    \end{equation*}
    For any $\widetilde{\Tt}_s\in\mathcal{E}$, we have for some $\lambda\in[0,1]$:
    \begin{equation*}
        \N{\widetilde{\Tt}_s}_{L^2(0,T;L^2(\Omega_s))} = \N{\lambda\mathcal{F}(\widetilde{\Tt}_s)}_{L^2(0,T;L^2(\Omega_s))} \leq \N{\mathcal{F}(\widetilde{\Tt}_s)}_{L^2(0,T;H^1_{\Gamma_D}(\Omega_s))}\leq C,
    \end{equation*}
    where the last estimate follows from \eqref{eq:theta_estimates}, see Proposition \ref{prop:theta_solution_estimates}. As a result, $\mathcal{E}$ is bounded in $L^2(0,T;L^2(\Omega_s))$. \par
    With Lemma \ref{lemm:F_continuity}, we conclude from Schaefer's fixed-point theorem that $\mathcal{F}$ has a fixed point $\Tt_s\in L^2(0,T;L^2(\Omega_s))$. Thus, there exists a weak solution $(c,w,q,\Tt_f,\Tt_s)$ of the fully coupled system \eqref{eq:solution_spaces}--\eqref{eq:coefficients}, which is given by Propositions \ref{prop:c_solution_estimates}, \ref{prop:w_solution_estimates}, and \ref{prop:theta_solution_estimates}.
    \qed
\end{proo}

\subsection{Uniqueness of solutions for the fully coupled system} \label{subsec:uniqueness}

In this section, we use the notation $\delta f=f^1-f^2$ for the difference of two functions $f^1$ and $f^2$. We will need the following Lipschitz estimates:

\begin{lemm} \label{lemm:Lipschitz_estimates}
    Let $\widetilde{\Tt}_s^1,\widetilde{\Tt}_s^2\in L^2(0,T;L^2(\Omega_s))$ and denote by
    \begin{equation*}
        (c^1,w^1,q^1,\Tt^1_f,\Tt^1_s)\qquad\text{and}\qquad(c^2,w^2,q^2,\Tt^2_f,\Tt^2_s)
    \end{equation*}
    the corresponding solutions of the linearized problem given by Propositions \ref{prop:c_solution_estimates}, \ref{prop:w_solution_estimates}, and \ref{prop:theta_solution_estimates}. Moreover, we write $(R^1,S_1,v^1_b,F_1,J_1,A_1)$ and $(R^2,S_2,v^2_b,F_2,J_2,A_2)$ for the corresponding quantities. Then there exists a constant $C>0$ such that for all $t\in (0,T]$, we have the following estimates:
    \begin{align}
        \N{\delta R}_{C_1^2([0,t]\times [0,L])} + \N{\delta S}_{C_1^2([0,t]\times \ov{\Omega})} + \N{\delta w}_{C([0,t];H^1_{\Sigma}(\Omega_f))} &\leq C\N{\delta\widetilde{\Tt}_s}_{L^2((0,t)\times\Omega_s)},\label{eq:solution_Lipschitz_estimates} \\
        \N{F}_{C^1([0,t]\times\ov{\Omega})^{3\times 3}} + \N{\delta J}_{C^1([0,t]\times\ov{\Omega})} + \N{\delta A}_{C^1([0,t]\times\ov{\Omega})^{3\times 3}} &\leq C\N{\delta\widetilde{\Tt}_s}_{L^2((0,t)\times\Omega_s)}, \label{eq:coefficients_Lipschitz_estimates} \\
        \N{\delta v_b}_{C([0,t]\times\ov{\Omega})^3} + \N{\G\delta v_b}_{C([0,t]\times\ov{\Omega})^{3\times 3}} &\leq C\N{\delta\widetilde{\Tt}_s}_{L^2((0,t)\times\Omega_s)}. \label{eq:v_b_Lipschitz_estimates}
    \end{align}
\end{lemm}

\begin{proo}
    The estimates can be derived with the same arguments as in the proof of Lemma \ref{lemm:F_continuity}.
    \qed
\end{proo}

\begin{proo}[Uniqueness for the fully coupled system \eqref{eq:solution_spaces}--\eqref{eq:coefficients}]
    Let
    \begin{equation*}
        (c^1,w^1,q^1,\Tt^1_f,\Tt^1_s)\qquad\text{and}\qquad(c^2,w^2,q^2,\Tt^2_f,\Tt^2_s)
    \end{equation*}
    be two weak solutions to the fully coupled system and denote by
    \begin{equation*}
        (R^1,S_1,v^1_b,F_1,J_1,A_1)\qquad\text{and}\qquad (R^2,S_2,v^2_b,F_2,J_2,A_2)
    \end{equation*}
    the corresponding quantities. We subtract \eqref{eq:transport_weak_ref} for each of the solutions, test the equation with $(\delta\Tt_f,\delta\Tt_s)$, and integrate over time from $0$ to some $t\in(0,T)$ to obtain:
    \begin{align*}
        &\int_0^t\langle\pa_t(J_1\Tt^1_f-J_2\Tt^2_f),\delta\Tt_f\rangle_{H^1(\Omega_f)',H^1(\Omega_f)}\,d\tau + \int_0^t\langle\pa_t(J_1\Tt^1_s-J_2\Tt^2_s),\delta\Tt_s\rangle_{H^1_{\Gamma_D}(\Omega_s)',H^1_{\Gamma_D}(\Omega_s)}\,d\tau \\
        & + \int_0^t\int_{\Omega_f}K_f^{F_2}\G\delta\Tt_f\cdot\G\delta\Tt_f\,dx\,d\tau + \int_0^t\int_{\Omega_s}K_s^{F_2}\G\delta\Tt_s\cdot\G\delta\Tt_s\,dx\,d\tau \\
        &= - \int_0^t\int_{\Omega_f}\delta K_f^{F}\G\Tt^1_f\cdot\G\delta\Tt_f\,dx\,d\tau - \int_0^t\int_{\Omega_s}\delta K_s^{F}\G\Tt^1_s\cdot\G\delta\Tt_s\,dx\,d\tau \\
        &\quad\underbrace{- \int_0^t \int_\Sigma \alpha|\delta\Tt_f-\delta\Tt_s|^2 J_1|F_1^{-T}n|\,d\sigma\,d\tau}_{\leq 0} \\
        &\quad \underbrace{- \int_0^t \int_\Sigma \alpha(\Tt^2_f-\Tt^2_s)(\delta\Tt_f-\delta\Tt_s)(J_1|F_1^{-T}n| - J_2|F_2^{-T}n|)\,d\sigma\,d\tau}_{=: I_1} \\
        &\quad \underbrace{- \int_0^t\int_{\Gamma_f} |\delta\Tt_f|^2(w^1\cdot n)^+ J_1\,d\sigma\,d\tau}_{\leq 0} \underbrace{- \int_0^t\int_{\Gamma_f} \Tt^2_f\delta\Tt_f((w^1\cdot n)^+ J_1 - (w^2\cdot n)^+ J_2)\,d\sigma\,d\tau}_{=:I_2} \\
        &\quad \underbrace{- \int_0^t\int_{\Gamma_f} f_{in}\delta\Tt_f((w^1\cdot n)^- J_1 - (w^2\cdot n)^- J_2)\,d\sigma\,d\tau}_{=:I_3} \\
        &\quad \underbrace{+ \int_0^t\int_{\Omega_f} \Tt^2_f (A_1w^1 - A_2w^2)\cdot\G\delta\Tt_f\,dx\,d\tau}_{=:I_4} \underbrace{+ \int_0^t\int_{\Omega_f} \delta\Tt_f A_1w^1\cdot\G\delta\Tt_f\,dx\,d\tau}_{=:I_5}.
    \end{align*}
    The first six terms have been treated, e.g., in the proof of \cite[Theorem 3]{gahn_homogenization_2023}. Note that Propositions \ref{prop:S_solution_estimates}, \ref{prop:w_solution_estimates}, and \ref{prop:theta_solution_estimates} and Lemma \ref{lemm:Lipschitz_estimates} guarantee that the same estimates as in the proof of \cite[Theorem 3]{gahn_homogenization_2023} also hold in our setting. From the last inequality therein, we obtain after absorbing the norm of the gradient by the left hand side: 
    \begin{align*}
        &c_1\left(\N{\delta\Tt_f(t)}^2_{L^2(\Omega_f)} + \N{\delta\Tt_s(t)}^2_{L^2(\Omega_s)}\right) + c_2\left(\N{\G\delta\Tt_f}^2_{L^2((0,t)\times\Omega_f)} + \N{\G\delta\Tt_s}^2_{L^2((0,t)\times\Omega_s)}\right) \\
        &\leq C \left(\N{\delta\Tt_f}^2_{L^2((0,t)\times\Omega_f)} + \N{\delta\Tt_s}^2_{L^2((0,t)\times\Omega_s)}\right) + \sum_{i=1}^5 I_i,
    \end{align*}
    where $c_1,c_2,C>0$ do not depend on $t$. \par
    For $I_5$, we use Hölder's inequality, an interpolation inequality, continuous embedding, and twice Young's inequality to obtain, for any $\nu>0$:
    \begin{align*}
        &\int_0^t\int_{\Omega_f} \delta\Tt_f A_1w^1\cdot\G\delta\Tt_f\,dx\,d\tau \\
        &\leq\widetilde{C}\int_0^t\N{\delta\Tt_f}^{\frac{1}{2}}_{H^1(\Omega_f)}\N{\delta\Tt_f}^{\frac{1}{2}}_{L^2(\Omega_f)}\N{A_1}_{C(\ov{\Omega_f})}\N{w^1}_{H^1(\Omega_f)}\N{\G\delta\Tt_f}_{L^2(\Omega_f)}\,d\tau \\
        &\leq \widetilde{C}\N{A_1}_{C([0,t]\times\ov{\Omega_f})}\N{w^1}_{C([0,t];H^1(\Omega_f))}\left(C_\nu \N{\delta\Tt_f}^2_{L^2((0,t)\times\Omega_f)} + \nu\N{\G\delta\Tt_f}^2_{L^2((0,t)\times\Omega_f)}\right),
    \end{align*}
    where $\widetilde{C}>0$ arises just from interpolation and embedding. \par
    To the terms $I_1$, $I_2$, $I_3$, and $I_4$, we apply the same estimates as in \eqref{eq:addiff_a_k_convergence} and \eqref{eq:addiff_b_k_convergence}, where for the boundary terms, that is $I_1$, $I_2$, and $I_3$, we use the following scaled trace inequality: For any $\ee>0$, there exists $C_\ee>0$ such that
    \begin{equation}
        \N{f}_{L^4(\pa\Omega_{f/s})}\leq C_\ee\N{f}_{L^2(\Omega_{f/s})} + \ee\N{\G f}_{L^2(\Omega_{f/s})}, \label{eq:scaled_trace}
    \end{equation}
    for all $f\in H^1(\Omega_{f/s})$, see, e.g., \cite[Exercise II.4.1]{galdi_introduction_2011}. Note that `$f/s$' indicates an equation for both `$f$' and `$s$'. In addition, we use the uniform upper bounds provided by Propositions \ref{prop:S_solution_estimates}, \ref{prop:w_solution_estimates}, and \ref{prop:theta_solution_estimates} as well as the Lipschitz estimates from Lemma \ref{lemm:Lipschitz_estimates} to obtain:
    \begin{align*}
        &c_1\left(\N{\delta\Tt_f(t)}^2_{L^2(\Omega_f)} + \N{\delta\Tt_s(t)}^2_{L^2(\Omega_s)}\right) + c_2\left(\N{\G\delta\Tt_f}^2_{L^2((0,t)\times\Omega_f)} + \N{\G\delta\Tt_s}^2_{L^2((0,t)\times\Omega_s)}\right) \\
        &\leq C \left(\N{\delta\Tt_f}^2_{L^2((0,t)\times\Omega_f)} + \N{\delta\Tt_s}^2_{L^2((0,t)\times\Omega_s)}\right) \\
        &\quad + C\N{\delta\Tt_s}_{L^2((0,t)\times\Omega_s)}\left(C_\ee\N{\delta\Tt_f}_{L^2((0,t)\times\Omega_f)} + \ee\N{\G\delta\Tt_f}_{L^2((0,t)\times\Omega_f)}\right) \\
        &\quad + C\N{\delta\Tt_s}_{L^2((0,t)\times\Omega_s)}\left(C_\ee\N{\delta\Tt_s}_{L^2((0,t)\times\Omega_s)} + \ee\N{\G\delta\Tt_s}_{L^2((0,t)\times\Omega_s)}\right) \\
        &\quad + C\N{\delta\Tt_s}_{L^2((0,t)\times\Omega_s)}\N{\G\delta\Tt_f}_{L^2((0,t)\times\Omega_f)} + C\left(C_\nu \N{\delta\Tt_f}^2_{L^2((0,t)\times\Omega_f)} + \nu\N{\G\delta\Tt_f}^2_{L^2((0,t)\times\Omega_f)}\right),
    \end{align*}
    where $C>0$ does not depend on $t$. Using again Young's inequality in each term on the right hand side but the first and the last, we can absorb all terms with gradients by the left hand side and obtain, with some new $C>0$, still not depending on $t$:
    \begin{equation*}
        \N{\delta\Tt_f(t)}^2_{L^2(\Omega_f)} + \N{\delta\Tt_s(t)}^2_{L^2(\Omega_s)} \leq C \left(\N{\delta\Tt_f}^2_{L^2((0,t)\times\Omega_f)} + \N{\delta\Tt_s}^2_{L^2((0,t)\times\Omega_s)}\right).
    \end{equation*}
    By Gronwall's inequality, we get $\N{\delta\Tt_f(t)}^2_{L^2(\Omega_f)} + \N{\delta\Tt_s(t)}^2_{L^2(\Omega_s)} = 0$ for a.e. $t\in(0,T)$, in particular $\Tt^1_s=\Tt^2_s$ a.e. in $(0,T)\times\Omega_s$. From the uniqueness of the solutions to all subproblems, see Propositions \ref{prop:c_solution_estimates}, \ref{prop:w_solution_estimates}, and \ref{prop:theta_solution_estimates}, it follows that the weak solution to the fully coupled system is unique.
    \qed
\end{proo}